\documentclass[reqno]{amsart} 
\usepackage[colorlinks]{hyperref}
\usepackage{amssymb}
\usepackage{latexsym}
\usepackage{amstext} 
\usepackage{color}
\usepackage{amsthm}  
\usepackage{amsmath}
\begin{document}
\newtheorem{theorem}{Theorem}    
\newtheorem{conjecture}{Conjecture}
\newtheorem{question}{Question}
\newtheorem{proposition}[theorem]{Proposition}
\newtheorem{lemma}[theorem]{Lemma}
\newtheorem{corollary}[theorem]{Corollary}
\newtheorem{example}[theorem]{Example}
\newtheorem{remark}[theorem]{Remark}
\newtheorem{definition}[theorem]{Definition}
\setlength{\parskip}{.15in} 
\setlength{\parindent}{0in} 

\title[\hfil  Abstract factorials]{Abstract factorials \\ } 
\author[Angelo B. Mingarelli \hfilneg]{Angelo B. Mingarelli}   
\address{School of Mathematics and Statistics\\ 
Carleton University, Ottawa, Ontario, Canada, K1S\, 5B6}
\email[A. B. Mingarelli]{amingare@math.carleton.ca}

\date{July 10, 2012 - Dedicated to the memory of my teacher, Professor Hans Heilbronn.}
\thanks{This research is partially supported by an NSERC Canada Discovery Grant}
\begin{abstract}
\noindent{A commutative semigroup of abstract factorials is defined in the context of the ring of integers. We study such factorials for their own sake, whether they are or are not connected to sets of integers. Given a subset $ X \subseteq\mathbb{Z^+}$ we construct a ``factorial set" with which one may define a multitude of abstract factorials on $X$. We study the possible equality of consecutive factorials, a dichotomy involving the limit superior of the ratios of consecutive factorials and we provide many examples outlining the applications of the ensuing theory; examples dealing with prime numbers, Fibonacci numbers, and highly composite numbers among other sets of integers. One of our results states that given any abstract factorial the series of reciprocals of its factorials always converges to an irrational number. Thus, for example, for any positive integer k the series of the reciprocals of the $k$-th powers of the cumulative product of the divisors of the numbers from 1 to n is irrational.}
\end{abstract}

\maketitle

\section{Introduction}

The study of generalized factorials has been a subject of interest during the past century culminating in the appearance of many treatises that fostered widespread applications. An introduction and review of this area has already been given in the recent paper \cite{mb1}  and so we shall not delve into the historical matter any further unless our applications require it so, referring the reader to \cite{mb1} for more information.

In this paper we restrict ourselves to the ring of integers only. Abstract factorials are defined as maps $!_{_a} : \mathbb{N} \to \mathbb{Z}^+$ satisfying very few conditions, conditions that are verified by apparently all existing notions of a generalized factorial in this context. 

As a consequence of the results herein we obtain, among other such results, the irrationality of the following numbers and classes of numbers, where $b, q, k \in \mathbb{Z^+}$ are arbitrary,
$$\sum_{n=1}^\infty \dfrac{1}{\prod_{i=1}^n {{p_i}^{k\,\lfloor n/i \rfloor}}}, \quad\sum_{n=1}^\infty \dfrac{1}{n!\,q^{\sum_{j=1}^n d(j)}}, \quad \sum_{n=1}^\infty \dfrac{1}{n!\,q^{\sum_{j=1}^n \sigma_k (j)}},$$
and 
$$\sum_{n=1}^\infty \frac{b^{nk}}{(bn)!^k},\quad  \sum_{n=1}^\infty \frac{1}{n!\,\mathcal{F}(n)^k}, \quad \sum_{n=1}^\infty \dfrac{1}{q^{\sum_{j=1}^n d(j)}\alpha(n)^k}\quad \sum_{n=1}^\infty \frac{1}{p_n!},$$
where $\alpha(n) = \prod_{i=1}^n i^{\lfloor n/i \rfloor} $ is the cumulative product of all the divisors from $1$ to $n$, $\mathcal{F}(n)$ is the product of the first $n$ Fibonacci numbers, and $p_n$ is the $n$-th prime.  Furthermore, there is a sequence of {\it highly composite numbers} $h_n$ such that 
$$\sum_{n=1}^\infty \dfrac{1}{\prod_{i=1}^n {{h_i}^{\lfloor n/i \rfloor}}}$$ is irrational. In addition, if $f:\mathbb{Z^+}\times \mathbb{Z^+}\to \mathbb{Z^+}$ satisfies a concavity condition in its first variable, i.e., for any $x,y,q$, we have $f(x+y,q) \geq f(x,q)+f(y,q)$, then, for any $q \in \mathbb{Z^+}$, 
$$\sum_{n=1}^\infty \frac{1}{q^{f(n,q)}\,n!}$$
is also irrational. We also show that we may choose $f(n, q) = {\binom{n+q-1}{q}}$ in the previous result.

An abstract factorial will be denoted simply by the notation $n!_{_a}$, the usual factorial function being denoted by $n!$. Other unspecified abstract factorials will be indexed numerically (e.g., $n!_{_1}, n!_{_2}, \ldots$.)

We always assume that $X$ is a non-empty set of non-zero integers. For the sake of simplicity assume that $X \subseteq \mathbb{Z^+}$, although this is, strictly speaking, not necessary as the constructions will show. Using the elements of $X$ we construct a new set (generally not unique) dubbed a {\it factorial set} of $X$.  We will see that any factorial set of $X$ may be used out to construct infinitely many abstract factorials (see Section~\ref{asso}). This construction of generalized (abstract) factorials of $X$ should be seen as complementary to that of Bhargava,  \cite{mb1}.  Furthermore, there is enough structure in the definition of these abstract factorials that, as a collective, they form a  semigroup under ordinary multiplication.

By their very nature abstract factorials should go hand-in-hand with binomial coefficients (Definition~\ref{def2}).  For example,  Knuth and Wilf  \cite{kw} define a generalized binomial coefficient by first starting with a positive integer sequence $\mathcal{C} = \{C_n\}$ and then defining the binomial coefficient as a ``falling" chain type product $$\binom{n+m}{m}_{\mathcal{C}}=C_{m+n}C_{m+n-1}\ldots C_{m+1}/C_nC_{n-1}\ldots C_1.$$ In this case, the quantity $n!_{_a}=n! C_1C_2\cdots C_n$ always defines an abstract factorial according to our definition. 

In Section~{\ref{s1}} we give the definition of an abstract factorial (Definition~\ref{def2}), give their representation (Proposition~\ref{binco})  and show that, generally, consecutive equal factorials may occur. In fact, strings of three or more consecutive equal factorials cannot occur (Lemma~\ref{10x}). Of special interest is the quantity defined by the ratio of consecutive factorials \eqref{bk} for which there exists a dichotomy, {\em i.e.,} there always holds either \eqref{eq100} or \eqref{eq101} (Lemma~\ref{lem000}). Cases of equality in both \eqref{eq100} or \eqref{eq101} are exhibited by specific examples (Proposition~\ref{excep} in the former case, and use of Bhargava's factorials for the set of primes \cite{mb1} in the latter case). 

Generally, given a set $X$ we find its factorial sets (Section~\ref{asso}). We then show that for any abstract factorial (whether or not it should arise from a set) the sum of the reciprocals of its generalized  factorials is always irrational (Section~\ref{s2}). An application of the semigroup property (Proposition~\ref{semi}) and the global irrationality result (from  Lemma~\ref{th2} and Lemma~\ref{th4}) implies that if $!_{1},$ $!_{2}$,\dots,$!_{k}$ is any collection of abstract factorials, $s_i \in \mathbb{N}$, $i=1,2,\dots,k$, not all equal to zero, then $$\sum^\infty_{n=0} \frac{1}{\prod^k_{j=1} n!^{s_j}_{j}} $$ is irrational (Theorem~\ref{thh}). As a consequence of our theory we also obtain the irrationality of the series of reciprocals of the generalized factorials (and their powers) for the set of primes in \cite{mb1}, (see Corollary~\ref{cor311}) and the other series displayed earlier.

In Section~\ref{ss5} we consider an inverse problem that may be stated thus: Given any abstract factorial $n!_{_a}$,  does there exist a set $X$ such that the sequence of generalized factorials $\{n!_{_a}\}_{n=0}^\infty$ coincides with one of the factorial sets of $X$? If there is such a set $X$, it will be called a {\it primitive} of the abstract factorial in question. It is noteworthy that such primitives are usually not unique.

In this direction we find that a primitive of the ordinary factorial function, $n!$, is given simply by the exponential of the Von Mangoldt function i.e., 
$ X = \{e^{\Lambda (m)} : m=1,2,\ldots\}$. In other words, the ordered set $$X=\{1, 1, 3, 1, 5, 1, 7, 1, 3, 1, 11, 1, 13, \ldots\}$$ whose $n$-th term is given by $b_n = e^{ \Lambda(n)}$ has a factorial set whose elements coincide with the sequence of factorials of the ordinary factorial function (Theorem~\ref{vonm}). We find (Theorem~\ref{mbpr}) that Bhargava's generalized factorial for the set of primes also has a primitive
$$X=\{2, 6, 1, 10, 1, 21, 1, 2, 1, 11, 1, 13,\ldots\}   $$
where every term here is the product of at most {\it two} primes. Still, it has a factorial function that agrees with the generalized factorial for the set of primes in \cite{mb1}. Thus, generally speaking, there are a number of ways in which one may associate a set with an abstract factorial and conversely.

In Section~\ref{s4} we give some applications of the foregoing theory.  We also introduce the notion of a self-factorial set, that is, basically a set whose elements are either the factorials of some abstract factorial, or are so when multiplied by $n!$. We find some abstract factorials of sets such as the positive integers, $X=\mathbb{Z^+}$, and show that one of its factorial sets is given by the set $\{n!_{_a}: n=0,1,2,\dots\}$ where (see Example~\ref{ex44}) $$n!_{_a} = \prod_{i=1}^n i^{i\,\lfloor n/i \rfloor}.$$
Furthermore, in Example~\ref{ex42} we show that one of the factorial sets of the set $\{1,q,q,q, \ldots\}$ where $q \in \mathbb{Z^+}$, $q\geq 2$, is given by the set $\{B_n\}$ where $$B_n = q^{\sum_{k=1}^n d(k)}$$ where $d(n)$ is the usual divisor function. 

Combining the preceding with the results of Section~\ref{s2} we also obtain that the series of reciprocals of the $k$-th powers ($k \geq 1$) of the cumulative product of all the divisors of the integers from $1$ to $n$, i.e., 
$$\sum_{n=1}^\infty 1/\prod_{i=1}^n i^{k\,\lfloor n/i \rfloor},$$ is irrational (see Example~\ref{ex44} and Remark~\ref{rem9}).

In the same spirit we show in Example~\ref{ex45} that for $q \in \mathbb{Z^+}$, $q\geq 2$, the set $\{q^n : n\in \mathbb{N}\}$ has a factorial set $\{B_n\}$ where $$B_n= q^{\sum_{k=1}^n \sigma(k)}$$ where $\sigma(n)$ is the sum of the divisors of $n$, a result that can be extended to the case of sets of integers of the form $\{q^{n^k}\}$ for given $k\geq 1$ (see Example~\ref{exnk}). Standard arithmetic functions abound in this context as can be gathered by considering the more general  situation $X=q\mathbb{Z}^+$, $q>0$. Here, one of the factorial sets of $X$ is given by numbers of the form $$B_n= q^{\sum_{k=1}^n d(k)}\, \prod_{i=1}^n i^{\lfloor n/i\rfloor},$$ where the product on the right is once again the cumulative product arithmetic function defined above (see Remark~\ref{rem9}).

Subsections~\ref{s5}-\ref{s62} are devoted to questions involving prime numbers in our set $X$, their (abstract) factorials, and the problem of determining whether a function arising from the ordinary factorial of the n-th prime number is, indeed, an abstract factorial. This latter question is, in fact, related to an unsolved problem of Hardy and Littlewood dealing with the convexity of the prime-counting function $\pi(x)$.

We show that the Hardy-Littlewood conjecture on the prime counting function, $\pi(x)$, i.e., that for all $x, y $ there holds
$$\pi(x+y) \leq \pi(x)+\pi(y),$$
implies that
$$p_n \geq p_k +p_{n-k-1}, \quad \quad 1 \leq k \leq n-1,$$
where $p_n$ is the $n$-th prime, and that this inequality in turn implies that the ``prime factorial function" $f : \mathbb{N}\to\mathbb{Z^+}$ defined by $f(0)=1$, $f(1)=1$ and $f(n) = p_{n-1}!$, $n \geq 2$, is an abstract factorial. Although said conjecture may be false according to some, it may be the case that the above inequality holds.

We recall the definition of a highly composite number (hcn): A number $n$ is said to be highly composite if $d(m) < d(n)$ whenever $m < n$, where $d$ is the usual divisor function. After proceeding to the calculation of a factorial set of the set of primes, we note that the first six numbers of this set are actually highly composite numbers and, in fact, we prove that these are the only ones (Proposition~\ref{fff}).

An application of the theory developed here allows us to derive that for every positive integer $k$, the series
$$\sum_{n=0}^{\infty} 1/(p_1^{\lfloor n/1\rfloor }p_2^{\lfloor n/2\rfloor}p_3^{\lfloor n/3\rfloor}\cdots p_n^{\lfloor n/n\rfloor})^k$$ is  irrational. 

In Subsection~\ref{s6} we show, in particular, that given any positive integer $m$ there is a highly composite number (hcn), $N$, such that $m! | N$. We then find factorial sets of the set of hcn and show that they are all self-factorial. Using this it is shown that there exists a sequence $\{h_n\}$ of hcn such that for any $k \in \mathbb{Z^+}$, 
$$\sum_{n=1}^\infty \dfrac{1}{\prod_{i=1}^n {{h_i}^{k\,\lfloor n/i \rfloor}}}$$ is irrational (Proposition~\ref{ffffff}). This is one of the few results dealing with the irrationality of series involving hcn.

We end the paper with a number of remarks. A brief note on the representations of abstract factorials in terms of possible solutions of the Stieltjes moment problem (Section~\ref{s8}) is given at the end. The idea here is to find integral representations of these abstract factorials akin to the usual representation of the ordinary factorial in terms of the Gamma function. Initial results in this direction indicate that various classes of abstract factorials admit unique integral representations as solutions of a moment problem. We note that not all abstract factorials admit such representations. We also produce a few more simple irrationality criteria based on the results herein, such as the ones given at the outset, and give a simple proof that the sum of the reciprocals of the factorials of the primes is irrational.

\section{Preliminaries}\label{s1}

In the sequel the symbols $X, I$ will always stand for non-empty subsets of $\mathbb{Z}$, not containing $0$, either may be finite or infinite, whose elements are not necessarily distinct (e.g., thus the set $X=\{1,q,q,q,\ldots\}$ is considered an infinite set). When the context is clear we will occasionally use the words sequence and sets interchangeably.
\begin{definition} \label{def2} An abstract (or generalized) factorial is a function $!_{_{a}}: \mathbb{N} \to \mathbb{Z^+}$ that satisfies the following conditions:
\begin{enumerate}
\item $0!_{_{a}}=1$,
\item For every non-negative integers $n, k$, $0 \leq k \leq n$ the generalized binomial coefficients $$\binom{n}{k}_{_{a}} := \frac{n!_{_{a}}}{k!_{_{a}}(n-k)!_{_{a}}} \in \mathbb{Z^+},$$
\item For every positive integer $n$, $n!$ divides $n!_{_{a}}$.
\end{enumerate}
\end{definition}

\begin{remark}\label{nondec} {\rm Since, by hypothesis (2) above, $\binom{n+1}{n}_{_{a}} \in \mathbb{Z^+}$ for every $n \in \mathbb{N}$ the sequence of abstract factorials $n!_{_{a}}$ is non-decreasing. }
\end{remark}
Another simple consequence of the definition is,
\begin{proposition} \label{semi} The collection of all abstract factorials forms a commutative semigroup under ordinary multiplication.
\end{proposition}

{\bf Terminology:}  In the sequel an abstract factorial function will be called simply a {\it factorial function} or an {\it abstract factorial} and its values will be referred to simply as its {\it factorials} (or generalized factorials for emphasis), unless otherwise specified.

Of course the ordinary factorial function $n!$ is an abstract factorial as is the function defined by setting $n!_{_a}:=2^{n(n+1)/2}n!$. The factorial function defined in \cite{mb1}, for arbitrary sets $X$ is also a factorial function (see Example~\ref{prop1}). In addition, if $\mathcal{C} = \{C_n\}$ is a positive integer sequence and we assume as in \cite{kw} that the binomial coefficient $$\binom{n+m}{m}_{\mathcal{C}}=\frac{C_{m+n}C_{m+n-1}\ldots C_{m+1}}{C_nC_{n-1}\ldots C_1}$$ is a positive integer for every $n, m\in \mathbb{N}$, then there is an associated abstract factorial $!_{_a}$ with these as binomial coefficients that is, the one defined by setting $0!_{_a}=1$ and $$n!_{_a} = C_1C_2\cdots C_n$$ provided $n! | C_1C_2\cdots C_n$ for every $n \in \mathbb{Z^+}$. On the other hand, if $n!$ does not divide $C_1C_2\cdots C_n$ for every $n$ we can still define another abstract factorial by writing $$n!_{_a} = n!\, C_1C_2\cdots C_n.$$ Its binomial coefficients are now of the form 
$$\binom{n+m}{m}_{_{a}} = \binom{n+m}{m}_{\mathcal{C}}\binom{n+m}{m}$$
where the last binomial coefficient is the usual one. These generalized or abstract binomial coefficients are necessarily integers because of the tacit assumption made in \cite{kw} on the binomial coefficients appearing in the middle of the previous display. All of our results below apply in particular to either one of the two preceding factorial functions.

Another consequence of Definition~\ref{def2} is the following,

\begin{proposition}\label{binco}
Let $!_{_a}$ be an abstract factorial. Then there is a positive integer sequence $h_n$ with $h_0=1$ and such that for each $n \in \mathbb{N}$,
\begin{equation}\label{bin1}
h_k h_{n-k} \bigg |\, h_n \binom{n}{k},\quad k=0,1,2\dots,n.
\end{equation}
Conversely, if there is a sequence of positive integers $h_n$ satisfying \eqref{bin1} and $h_0=1$, then the function $!_{_a}: \mathbb{N}\to \mathbb{Z^+}$ defined by $n!_{_a}= n!h_n$ is an abstract factorial.
\end{proposition}

\begin{corollary} Let $h_n \in \mathbb{Z^+}$ be such that $h_0=1$, $h_k h_{n-k}  | h_n$, for all $k=0,1,\dots,n,$ and for every $n \in \mathbb{N}$. Then $n!_{_a}= n!h_n$ is an abstract factorial.
\end{corollary}

In Section~\ref{asso} below we consider those abstract factorials induced by those sequences $h_n$ such that  $h_k h_{n-k}  | h_n$, for all $k=0,1,\dots,n,$ and for every $n \in \mathbb{Z^+}$. Such sequences form the basis for the notion of a ``factorial set" as we see below.

Observe that $h_n$ is a constant sequence satisfying \eqref{bin1} if and only if $h_n =1$ for all $n$, that is, if and only if the abstract factorial reduces to the ordinary factorial. Note that {\it distinct abstract factorials functions may have the same set of binomial coefficients}; for example, if $b \in \mathbb{Z^+}$ and $n!_{_a} = n! b^n$, for every $n$, then the binomial coefficients of this factorial function and the usual factorial function are identical.

The reason for this lies in the easily verifiable identity
$$ n!_{_a} = 1!_{_a}^n\,\, \prod_{m=1}^n \binom{m}{m-1}_{_a},$$
valid for any abstract factorial. Thus, it is the value of $1!_{_a}$ that determines whether or not an abstract factorial is determined uniquely by a knowledge of its binomial coefficients.

One of the curiosities of abstract factorials lies in the possible existence of {\it equal consecutive factorials}. 

\begin{definition}\label{cons}
Let $!_{_a}$ be an abstract factorial. By a pair of equal consecutive factorials we mean a pair of consecutive factorials such that, for some $k\geq 2$, $k!_{_a}=(k+1)!_{_a}$. 
\end{definition}
\begin{remark} {\rm Definition~\ref{cons} is not vacuous as we do not tacitly assume that the factorials form a strictly increasing sequence (cf., Example~\ref{ex10} and Proposition~\ref{excep} below). In addition, given an abstract factorial it is impossible for  $1!_{_a}=2!_{_a}$. (This is because $\binom{2}{1}_{_a}$ must be an integer, which of course can never occur since $2!_{_a}$ must be at least $2$.)
}
\end{remark}

Such equal consecutive factorials, when they exist, are connected to the properties of ratios of {\it nearby} factorials. We adopt the following notation for ease of exposition: For a given integer $k$ and for a given factorial function $!_{_a}$, we write 
\begin{equation}\label{bk}
r_k = \frac{(k+1)!_{_a}}{k!_{_a}}.
\end{equation} 
Since generalized binomial coefficients are integers by  Definition~\ref{def2}, $r_k$ is an integer for every $k =0,1,2,\ldots$. The next result shows that strings of {\rm three or more} equal consecutive factorials cannot occur.
\begin{lemma}\label{10x} There is no abstract factorial with three consecutive equal factorials.
\end{lemma}
Knowing that consecutive equal factorials must occur in pairs if they exist at all we get,
\begin{lemma}\label{11x} Given an abstract factorial $!_{_a}$, let $2!_{_a} \neq 2$. If $r_k=1$ for some $k\geq 2$, then $r_{k-1}\geq 3$.
\end{lemma}
The next result gives a limit to the asymptotics of sequences of ratios of consecutive factorials defined by the reciprocals of the $r_k$. These ratios do not necessarily tend to zero as one may expect (as in the case of the ordinary factorial), but may have subsequences approaching non-zero limits!
\begin{lemma} \label{lem000} For any abstract factorial, either 
\begin{equation}\label{eq100}
\limsup_{k \to \infty}  \frac{1}{r_k}  = 1,
\end{equation}
or 
\begin{equation}\label{eq101}
\limsup_{k \to \infty}  \frac{1}{r_k}  \leq  1/2,
\end{equation}
the upper bound in \eqref{eq101} being sharp, equality being attained in the case of Bhargava's factorial for the set of primes (see the proof of Corollary~\ref{cor311}).
\end{lemma}
\begin{definition}\label{exc}{\rm An abstract factorial whose factorials satisfy \eqref{eq100} will be called {\bf exceptional.}}
\end{definition}
{\bf Note:} Using the generalized binomial coefficients $\binom{n+1}{n}_{_{a}}$ it is easy to see that a necessary condition for the existence of such exceptional factorial functions is that $1!_{_a} =1$. The question of their existence comes next.

\begin{proposition}\label{excep} The function $!_{_a} :\mathbb{N}\to \mathbb{Z^+}$ defined by $0!_{_a}=1$, $1!_{_a}=1$ and inductively by setting $ (n+1)!_{_a} = n!_{_a} $ whenever $n$ is of the form $n=3m-1$ for some $m \geq 1$, and 
\[n!_{_a} = \left \{ \begin{array}{ll}
		 n!\,(n+1)! \prod_{j=1}^{n-1}{(n - j)!_{_a}}^2, &\mbox{if  $n$ is of the form $n=3m-1$}, \\
		 n!\,\prod_{j=1}^{n-1}{(n - j)!_{_a}}^2, &\mbox{if $n$ is of the form $n=3m+1$}.  \\
		\end{array}
	\right.
\]
 is an exceptional factorial function.
\end{proposition}

\begin{remark} {\rm The construction in Proposition~\ref{excep} may be generalized simply by varying the exponent outside the finite product from 2 to any arbitrary integer greater than two. There then results an infinite family of such exceptional factorials. The quantity defined by $\prod_{j=1}^{n-1}(n - j)!_{_a}$, may be thought of as an abstract generalization of the {\it super factorial} (see \cite{njs}, A000178)}.
\end{remark}
\begin{example} The first few terms of the exceptional factorial defined in Proposition~\ref{excep} are given by $1!_{_a}=1$, $2!_{_a}=3!_{_a}=12$, $4!_{_a}=497664$, $5!_{_a} = 6!_{_a} = 443722221348087398400$, etc.
\end{example}
Since the preceding results are valid for abstract factorials they include, in particular, the recent factorial function considered in \cite{mb1}, and we summarize its construction for completeness. Let $X \subseteq {\bf Z}$ be a finite or  infinite set of integers. Following \cite{mb1}, we define the notion of a $p$-ordering of $X$ and use it to define the generalized factorials of the set $X$ inductively. By definition $0!_X=1$. For $p$ a prime, we fix an element $a_0 \in X$ and, for $k \geq 1$, we select $a_k$ such that the highest power of $p$ dividing $\prod_{i=0}^{k-1}{(a_k-a_i)}$ is minimized. The resulting sequence of $a_i$ is called a $p$-ordering of $X$. As one can gather from the definition, such $p$-orderings are not unique, as one can vary $a_0$. Associated with such a $p$-ordering of $X$ we define an associated $p$-sequence $\{\nu_k(X,p)\}_{k=1}^{\infty}$ by $$\nu_k(X,p) = w_p(\prod_{i=0}^{k-1}{(a_k-a_i)}),$$ where $w_p(a)$ is, by definition, the highest power of $p$ dividing $a$ (e.g., $w_3(162)=81$). It is then shown  that although the $p$-ordering is not unique the associated $p$-sequence is independent of the $p$-ordering being used. Since this quantity is an invariant, one can use this to define generalized factorials of $X$ by setting \begin{equation}\label{mba} k!_X = \prod_{p}\nu_k(X,p),\end{equation} where the (necessarily finite) product extends over all primes $p$. 
\begin{example} \label{prop1} Bhargava's factorial function \eqref{mba} is an abstract factorial.
\end{example}
Hypothesis 1 of Definition~\ref{def2} is clear by definition of the factorial in question. Hypothesis 2 of Definition~\ref{def2} follows by the results in  \cite{mb1}. 

As we mentioned above, the question of the possible existence of equal consecutive generalized factorials is of interest. We show herewith that such examples exist for abstract factorials over the ring of integers.

\begin{example}\label{ex10} {\rm There exist sets $X$ with consecutive equal Bhargava factorials, $!_{_X}$. Perhaps the easiest example of such an occurrence lies in the set of generalized factorials of the set of cubes of the integers, $X=\{n^3 : n \in \mathbb{N}\}$, where one can show directly that $3!_X=4!_X  (=504)$. Actually, the first occurrence of this is for the finite subset $\{0, 1, 8, 27, 64, 125, 216, 343\}$.

Another such set of consecutive equal generalized factorials is given by the finite set of Fibonacci numbers $X=\{F_2,F_3,\ldots,F_{18}\}$, where one can show directly that $7!_{_a}=8!_{_a} (=443520)$. We point out that the calculation of factorials for finite sets as defined in \cite{mb1} is greatly simplified through the use of Crabbe's algorithm \cite{ac}. }
\end{example}
Inspired by the factorial representation of the base of the natural logarithms, one of the basic objects of study here is the series defined by the sum of the reciprocals of the factorials in question.
\begin{definition} For a given abstract factorial we define the constant $e_{_a}$ by the series of reciprocals of its factorials, i.e.,
\begin{equation}\label{eqqq1}
e_{_a} \equiv \sum_{r=0}^{\infty}{\frac{1}{{r!}_{_a}}}.
\end{equation}
\end{definition}
Note that the series appearing in \eqref{eqqq1} converges on account of Definition~\ref{def2}(3) and $1< e_{_a} \leq e$. 

\section{Factorial sets and their properties}\label{asso}
Besides creating abstract factorials using clever constructions, the easiest way to generate them is by means of integer sequences. As we referred to earlier it is shown in \cite{mb1} that on every subset $X\subseteq \mathbb{Z^+}$ one can define an abstract  factorial. We show below that there are other (non-unique) ways of generating abstract factorials (possibly infinitely many) out of a given set of positive integers. 

\subsection{The construction of a factorial set}

Given $I=\{b_1,b_2,\ldots\}$, $I \subset \mathbb{Z}$, ($b_i\neq 0$), with or without repetitions, we associate to it another set $X_I=\{B_0,B_1,$ $\ldots, B_n, \ldots \}$ of positive integers, termed simply a {\it factorial set of $I$}. In this case $I$ is called a {\it primitive} (set) of $X_I$. 

The elements of this factorial set $X_I $ are defined as follows: $B_0=1$ by definition, $B_1$ is (the absolute value of) an arbitrary but fixed element of $I$, say, $B_1=|b_1|$ (so that the resulting factorial set $X$ generally depends on the choice of $b_1$). Next, $B_2$ is the smallest (positive) number of the form ${b_1}^{\alpha_1}\,{b_2}^{\alpha_2}$ (where the $\alpha_i > 0$) such that ${B_1}^2 | B_2$. Hence $B_2=|{b_1}^{2}\,{b_2}|$. Next, $B_3$ is defined as the smallest (positive) number of the form ${b_1}^{\alpha_1}\,{b_2}^{\alpha_2}{b_3}^{\alpha_3}$ such that $B_1B_2 | B_3$. Thus, $B_3 = |{b_1}^{3}\,{b_2}\,{b_3}|$. Now, $B_4$ is defined as that smallest (positive) number of the form $\prod_{k=1}^4 {b_k}^{\alpha_k}$ such that $B_1B_3 | B_4$ and ${B_2}^2 | B_4$. This calculation gives us $B_4 = |{b_1}^{4}\,{b_2}^2\,{b_3}{b_4}|$. In general, we build up the elements $B_i$, $i=2,3,\ldots,n-1,$ inductively as per the preceding construction and define the element $B_n$ as that smallest (positive) number of the form $| \prod_{k=1}^n {b_k}^{\alpha_k}|$ such that $B_iB_j | B_n$  for every $i,j$, $0 \leq i \leq j \leq n$, and $i+j = n$. 

\begin{remark} {\rm Observe that permutations of the set $I$ may lead to different factorial sets, $X_I$, each one of which will be used to define a different abstract factorial (below). }
\end{remark}

It is helpful to think of the elements $B_n$ of a factorial set as defining a sequence of {\it generalized factorials}. In \cite{mb1}  one finds that the set of ordinary factorials arises from a general construction applied to the set of positive integers. For the analogue of this result see Section~\ref{ss5}.

The basic properties of any one of the factorial sets of a set of integers, all of which follow from the construction, can be summarized as follows.
\begin{remark}{\rm Let $I=\{b_i\} \subset \mathbb{Z}$  be any infinite subset of non-zero integers. For any fixed $b_m \in I$, consider the (permuted) set $I' = \{b_m, b_1, b_2, \ldots, b_{m-1}, b_{m+1}, \ldots\}$. Then the factorial set $X_{b_m}=\{B_1, B_2, \ldots, B_n, \ldots \}$ of $I'$ exists and for every $n > 1$ and for every $i,j \geq 0$, $i+j = n$, we have $B_iB_j | B_n$. In addition, if the elements of $I$ are all positive, then the $B_i$ are monotone.}
\end{remark}

Of course, factorial sets may be finite (e.g., if $X$ is finite) or infinite. The next result shows that factorial sets may be used to construct infinitely many abstract factorials.

\begin{theorem}\label{thh9} Let $I$ be an integer sequence, $X_I=\{B_n\}$ one of its factorial sets. Then, for each $k \geq 0$, the map $!_{_a}: \mathbb{N} \to \mathbb{Z^+}$ defined by ($0!_{_a}=1$)
$$ n!_{_a} = n! B_1B_2 \cdots B_{n+k},$$
is an abstract factorial on $I$.
\end{theorem}
Varying $k$ of course produces an infinite family of abstract factorials on $I$. The above construction of a factorial set leads to very specific sets of integers, sets whose elements we characterize next. (In the sequel, as usual,  ${\lfloor x \rfloor}$ is the greatest integer not exceeding $x$.) 
\begin{theorem}\label{th3} {\rm Given $I=\{b_i\} \subset \mathbb{Z^+}$, the terms 
\begin{equation}\label{bn} B_n = {b_1}^{\lfloor n \rfloor}{b_2}^{\lfloor n/2 \rfloor}{b_3}^{\lfloor n/3 \rfloor}
\cdots {b_n}^{\lfloor n/n \rfloor}
\end{equation}
characterizes one of its factorial sets, $X_{b_1}$. 
}
\end{theorem}
The next result leads to a structure theorem for generalized binomial coefficients corresponding to factorial functions induced by factorial sets.
\begin{proposition}\label{probi}
With $B_n$ defined as in \eqref{bn} we have, for every $n \in \mathbb{Z^+}$ and for every $k=0,1,\dots,n$,
\begin{equation}\label{bin01}
\frac{B_n}{B_k B_{n-k}} = \prod^n_{i=1} {b_i}^{\alpha_i},\quad \alpha_i =0\, \mbox{or}\,\, 1.
\end{equation}
\end{proposition}
With this the next result is clear.
\begin{corollary}\label{thebi} Let $n! | B_n$ for all $n\in \mathbb{Z^+}$. Then $n!_{_a} = B_n$ is an abstract factorial whose generalized binomial coefficients are of the form
\begin{equation}\label{bin02}
\binom{n}{k}_a = \prod^n_{i=1} {b_i}^{\alpha_i},\quad \alpha_i =0\, \mbox{or}\,\, 1.
\end{equation}
\end{corollary}

\section{Irrationality results}\label{s2}

We now state a few lemmas leading to a general irrationality result for sums of reciprocals of abstract factorials. First we note that given a positive integer sequence $b_n$ the series
\begin{equation}\label{nbn}\sum_{n=0}^\infty \dfrac{1}{n!b_n}\end{equation} may be either rational or not. Indeed, Erd\"{os} \cite{e75} pointed out that the series $$\sum_{n=0}^\infty \dfrac{1}{n!(n+2)} = 1.$$ The problem in this section consists in determining irrationality criteria for series of the form \eqref{nbn} using abstract factorials.
\begin{lemma}\label{th2} Let $!_{_a}$ be an abstract factorial whose factorials satisfy \eqref{eq101}. Then $e_{_a}$ is irrational.
\end{lemma}
\begin{remark} {\rm Although condition (3) in Definition~\ref{def2} (i.e., $n! | n!_{_a}$) of an abstract factorial appears to be very stringent, one cannot do without something like it; that is Lemma~\ref{th2} above is false for generalized factorials not satisfying this or some other similar property. For example, for $q > 1$ an integer, define the function $n!_{_a}=q^n$. It satisfies properties (1) and (2) of Definition~\ref{def2} but not (3). In this case it is easy to see that even though our function satisfies equation  \eqref{eq101}, $e_{_a}$ so defined is rational.}
\end{remark}
\begin{corollary}\label{cor311} Let $X$ be the set of prime numbers and $!_{_a}$ the factorial function \cite{mb1} of this set given by {\rm \cite{mb1} }
\begin{equation}\label{mbp}n!_{_a} = \prod_{p}^{} p^{\sum_{m=0}^{\infty}{[\frac{n-1}{p^m(p-1)}]}},
\end{equation}
where the (finite) product extends over all primes. Then $e_{_a}\approx 2.562760934$ is irrational.
\end{corollary}
The previous result holds because the generalized factorials of the set of primes satisfy \eqref{eq101} with equality. The next lemma covers the logical alternative exhibited by equation \eqref{eq100} in Lemma~\ref{lem000}.
\begin{lemma}\label{th4}  Let $!_{_a}$ be an abstract factorial whose factorials satisfy \eqref{eq100}. Then $e_{_a}$ is irrational. 
\end{lemma}
Combining the previous two lemmas we find the following theorem,
\begin{theorem}\label{thh} For any abstract factorial, $!_{_a}$ the number $e_{_a}$ is irrational. In fact, if $!_{a_1},$ $!_{a_2}$,\dots,$!_{a_k}$ is any collection of factorial functions, $s_i \in \mathbb{N}$, $i=1,2,\dots,k$, not all equal to zero, then $$\sum^\infty_{n=0} \frac{1}{\prod^k_{j=1} n!^{s_j}_{a_j}}$$
is irrational.
\end{theorem}

The alternating series counterpart of the previous theorem is next.
\begin{theorem}\label{thh2} Let $!_{a_1},$ $!_{a_2}$,\dots,$!_{a_k}$ be any collection of factorial functions, $s_i \in \mathbb{N}$, $i=1,2,\dots,k$, not all equal to zero. Then $$\sum^\infty_{n=0} \frac{(-1)^n}{\prod^k_{j=1} n!^{s_j}_{a_j}}$$
is irrational.
\end{theorem}
\begin{remark} {\rm This shows that the irrationality of the constants $e_{a}$ appears to be due more to the structure of the abstract factorial in question than an underlying theory about the base of the natural logarithms. }
\end{remark}
\begin{example}\label{ex44}{\rm  Let $n!_{_a} := (2n)!/2^n$, $n=0,1,2\ldots.$ Then this is an abstract factorial. An immediate application of Theorems~\ref{thh} and \ref{thh2} in the simplest case where $s=1$ gives that both quantities 
$$\sum_{n=0}^{\infty} \frac{2^n}{(2n)!} = \cosh \sqrt{2},$$ and 
$$\sum_{n=0}^{\infty} \frac{(-1)^n 2^n}{(2n)!} =  \cos \sqrt{2}$$ are irrational.

More generally, for $ b\in \mathbb{Z^+}$ the quantity $n!_{_b}=(bn)!/b^n$, $0!_{b}=1$,  defines an abstract factorial so that generally, $$ \sum_{n=0}^{\infty} \frac{b^n}{(bn)!}$$ is irrational as well along with its alternating series counterpart.}
\end{example}
\begin{example}{\rm Let $F_n$ denote the classical Fibonacci numbers defined by the recurrence relation $F_n = F_{n-1}+F_{n-2}$, $F_0=F_1=1$. The ``Fibonacci factorials" (\cite{njs}, id:A003266), denoted here by ${\mathcal{F}}(n)$ are defined by 
$${\mathcal{F}}(n) = \prod_{k=1}^n F_k.$$ Define a factorial function by setting ${\mathcal{F}}(0):=1$ and 
$$n!_{_a} := n! {\mathcal{F}}(n),\quad n=1,2,\ldots $$
In this case, the generalized binomial coefficients involve the {\it Fibonomial coefficients} ($={\mathcal{F}}(n)/{\mathcal{F}}(k){\mathcal{F}}(n-k)$) so that
$$\binom{n}{k}_{_a} = \binom{n}{k}\, \binom{n}{k}_{_F}$$ where the Fibonomial coefficients on the right (\cite{njs}, id:A010048), \cite{kw},  are integers for $k=0,1,\ldots,n$. Once again, an application of Theorem~\ref{thh} yields that for every $k \in \mathbb{Z^+}$,
 $$\sum_{n=0}^{\infty} \frac{1}{n!({\mathcal{F}}(n))^k} $$ 
is irrational.}
\end{example}
\begin{example}{\rm The (exceptional) abstract factorial of Proposition~\ref{excep} gives the rapidly growing (irrational) series of reciprocals of factorials:}
\begin{eqnarray*}
e_{_a} = 1+1+\frac{1}{12} + \frac{1}{12} + \frac{1}{497664} + \frac{1}{443722221348087398400} + \cdots \\ \\
\end{eqnarray*}
\end{example}

\section{An inverse problem}\label{ss5}

We recall that a set $X$ is called a primitive of an abstract factorial $!_{_a}$ if the sequence of its factorials, $\{n!_{_a}\}_{n=0}^\infty$ coincides with one of the factorial sets of $X$.  The question we ask here is: {\it When does an abstract factorial admit a primitive set?} Firstly, we give a simple necessary and sufficient condition for the existence of such a primitive set and, secondly, we give examples, the first of which shows that the ordinary factorial function has a primitive set whose elements are simply given by the exponential of the Von Mangoldt function.

\begin{lemma}\label{lns1}  A necessary and sufficient condition that a set $X=\{b_n\}$ be a primitive of the abstract factorial $n!_{_a}$ is that the quantity
\begin{equation}\label{eqlns1}
b_n = \frac{n!_{_a}}{\prod_{i=1}^{n-1} b_i^{\lfloor n/i \rfloor}    }
\end{equation}
defined recursively starting with $b_1=1!_{_a}$, be an integer for every $n > 1$.
\end{lemma}

\begin{remark}{\rm  It is not the case that \eqref{eqlns1} is always an integer even though the first three terms $b_1, b_2, b_3$ are necessarily so. The reader may note that the abstract factorial defined by  $n!_{_a} = (2n)!/2^n$ has no primitive since $b_4=14/3$. On the other hand, the determination of classes of abstract factorials that admit primitives is a fascinating problem.}
\end{remark}

We cite two examples of important factorials that do admit primitives.

\begin{theorem}\label{vonm} The ordinary factorial function has for a primitive (besides the set $\mathbb{Z^+}$) the set $X=\{b_n\}$ where $b_n = e^{\Lambda(n)}$, where $\Lambda(n)$ is the Von Mangoldt function.
\end{theorem}

Another example of an abstract factorial that admits a primitive (other than the original set it was intended for)  is the factorial function \cite{mb1} for the set of primes. In other words, there is a set $X$ different from the set of primes whose factorials (as defined herein) coincide with the abstract factorial 
\begin{equation}\label{eqmbpr}n!_{_a} = (n+1)!_X = \prod_{p}^{} \displaystyle p^{\sum_{m=0}^{\infty}{[\frac{n}{p^m(p-1)}]}},
\end{equation}
obtained in \cite{mb1} for the set of primes. (The factorial there is denoted  by $n!_X$).)

\begin{theorem}\label{mbpr} The abstract factorial defined in \eqref{eqmbpr} has for a primitive (besides the unordered set of prime numbers) the ordered set $X=\{b_i\}$ where here $b_1 = 2$, and the remaining $b_n$ are given recursively by \eqref{eqlns1} and explicitly as follows:
\[ b_i = \left \{ \begin{array}{ll}
		\quad\quad 1 ,	&	\mbox{if \ $i\neq p^m(p-1)$ for any prime $p$ and any $m\geq 0$}, \\
		\displaystyle {\prod_{p,\ \ i=p^m(p-1)} p} ,	&	\mbox{if \ $i= p^m(p-1)$ for some prime $p$ and $m \geq 0$,} \\
		\end{array}
	\right. 
\]
where the product extends over all primes $p$ such that $i$ has a representation in the form $i=p^m(p-1)$, for some $m \geq 0$.

\end{theorem}
The first few terms of the set $X$ in Theorem~\ref{mbpr} are given by 
\begin{eqnarray*}
X &=& \{2, 6, 1, 10, 1, 21, 1, 2, 1, 11, 1, 13, 1, 1, 1, 34, 1, 57, 1, 5, 1, 23,\\ & & 1, 1, 1, 1, 1, 29, 1, 31, 1, 2, 1, 1, 1, 37, 1, 1, 1,\ldots\}
\end{eqnarray*}
It follows from Theorem~\ref{mbpr} that if $i$ is odd then $b_i=1$, necessarily. It is tempting to conjecture that every term in $X$ is the product of at most two primes and this is, in fact, true. 

\begin{proposition}\label{dio} Let $n \in \mathbb{Z^+}$. Then there are at most two (2) representations of $n$ in the form $n=p^m(p-1)$ where $p$ is prime and $m \in \mathbb{N}$.
\end{proposition}

\section{Applications}\label{s4}

Before proceeding with some applications we require a few basic lemmas, the first of which, not seemingly well-known, is actually due to Hermite (\cite{her}, p.316) and rediscovered a few times since. e.g., Basoco (\cite{mab}, p.722, eq. (16).)
\begin{lemma}\label{lkp} For $k \geq 0$ an integer, let $\sigma_k(n)$ denote the sum of the $k$-th powers of the divisors of $n$, (where, $\sigma_0(n) = d(n)$). Then 
\begin{equation}\label{kp}
\sum_{i=1}^{n} \sigma_k(i) = \sum_{i=1}^{n} i^k\, \lfloor {n}/{i} \rfloor.
\end{equation}
\end{lemma}
{\bf Note:} The left-side of \eqref{kp} is the summatory function of $ \sigma_k(i)$ or $n$ times the average order of $ \sigma_k(i)$ over its range (\cite{hw}, Section 18.2). Furthermore, there is an interesting relationship between \eqref{kp} and the Riemann zeta function at the positive integers, that is,
$$\sum_{i=1}^{n} \sigma_k(i) = \dfrac{\zeta(k+1)}{k+1}n^{k+1} + O(n^k)$$
where the remainder terms are in terms of Ramanujan sums.
\begin{lemma}\label{lcp} Let $\alpha (n)$ denote the cumulative product of all the divisors of the numbers $1,2,\ldots, n$. Then 
\begin{equation}\label{ckp}
\alpha (n)  = \prod_{i=1}^{n} i^{\lfloor {n}/{i} \rfloor.}
\end{equation}
\end{lemma}
\begin{remark}{\rm It is also known that 
$$\alpha (n)  = \prod_{k=1}^n \, \lfloor  \frac{n}{k} \rfloor ! $$
(see \cite{njs}, id.A092143, Formula).}
\end{remark}
We now move on to examples where we describe explicitly some of the factorial sets of various basic integer sequences.
\begin{example}\label{ex42} {\rm The factorial set $X_I$ of the set $I$ of basically identical integers, $I=\{1,q,q,q,q, \ldots\}$ as per our construction where $q \geq 2$, and $B_1=q$, gives the factorial set
\begin{equation}
X_{I}=\{1,q, q^3, q^5, q^8, q^{10}, q^{14}, q^{16}, q^{20}, q^{23}, q^{27}, q^{29}, q^{35}, \ldots\}
\end{equation}
a set whose $n$-th term is $B_n = q^{a(n)}$, where $a(n) =\sum_{k=1}^n\,d(k)$ (by Theorem~\ref{th3} and Lemma~\ref{lkp}) and $d(k)$ is, as before, the number of divisors of $k$. The function defined by setting $n!_{_{a}} =n!\, B_n$  defines an abstract factorial. Here we see that equal consecutive factorials cannot occur by construction. In addition,  by Lemma~\ref{th2},
$$\sum_{n=1}^\infty \frac{1}{n!\,q^{\sum^n_{k=1} d(k)}} $$
is irrational.}
\end{example}
\begin{definition}
{\rm Let $I$ be an infinite subset of $\mathbb{Z^+}$ with a corresponding factorial set $X_I=\{B_n\}$. If $n! | B_n$ for every $n$, we say that this factorial set $X_I$ is a {\bf self-factorial set}.}
\end{definition}
The motivation for this terminology is that the function defined by setting $n!_{_{a}} = B_n$ is an abstract factorial. In other words, a self-factorial set may be thought of as an infinite integer sequence of consecutive generalized factorials (identical to the set itself, up to permutations of its elements). The next result is very useful when one wishes to iterate the construction of a factorial set {\it ad infinitum} (i.e., when finding the factorial set of a factorial set, etc.).
\begin{lemma}\label{asi}  If $I=\{b_n\}$ is a set with $n! | b_n$ for every $n$, then its factorial set $X_{b_1}$ is a self-factorial set. 
\end{lemma}
The same idea may be used to prove that
\begin{corollary}\label{cor46} The factorial set $X_{B_1}$ of a self-factorial set $X=\{B_n\}$ is a self-factorial set.
\end{corollary}
Next, we show that set $\mathbb{Z^+}$ has a factorial set with interesting properties. 
\begin{example}\label{ex44} {\rm We find a factorial set of the set  $X=\mathbb{Z^+}$ as per the preceding construction. Choosing $B_1=1$ we get the following set, 
\begin{equation}\label{div}
X_\mathbb{Z^+}=\{1,2,6,48,240,8640,60480,3870720,104509440,10450944000,\ldots\}
\end{equation}
a set which coincides (by Lemma~\ref{lcp} and Theorem~\ref{th3}) with the set of cumulative products of all the divisors of the numbers $1,2,\ldots,n$ (see Sloane \cite{njs}, id.A092143). Note that by construction $n! | B_n$ for every $n$. Hence, we can define an abstract factorial by setting $n!_{_{a}} = B_n$ to find that for this factorial function the set of factorials is given by the {\it set itself}, that is, this $X_\mathbb{Z^+}$ is self-factorial. In particular, equal consecutive factorials cannot occur by construction, and it follows from Lemma~\ref{th2} that the number defined by the sum of the reciprocals of these $B_n$, i.e., 
$$e_{_{a}}= 1+\sum_{n=1}^{\infty}\, 1/ \prod_{i=1}^{n} i^{\lfloor {n}/{i} \rfloor} = 1 + 1 + \frac{1}{2} + \frac{1}{6} + \frac{1}{48} + \frac{1}{240} + \ldots \approx 2.69179920$$
is irrational. Observe that infinitely many other integer sequences $I$ have the property that $n! | B_n$ for all $n$. Such sequences can thus be used to define abstract factorials. For example, if we consider the set of all $k$-th powers of the integers, $ I = \{n^k: n\in \mathbb{Z^+}\}, \, k\geq 2$, then another application of Lemma~\ref{lcp} shows that its factorial set $X_I$ (with $B_1=1$) is given by terms of the form $$B_n =  \prod_{i=1}^{n} i^{k\lfloor {n}/{i} \rfloor}.$$ In these cases we can always define an abstract factorial by writing $n!_{_{a}} = B_n$. Indeed, the semigroup property of abstract factorials (Proposition~\ref{semi}) implies that for each $k \in \mathbb{Z^+}$ the series of $k$-th powers of the reciprocals of this cumulative product,
$$\sum_{n=1}^\infty 1/\prod_{i=1}^n i^{k\,\lfloor n/i \rfloor} $$ is irrational.
\begin{remark}\label{rem9} {\rm The previous results are a special case of a more general result which states that the factorial set of the set $X=q\mathbb{Z^+}$, $q \in \mathbb{Z^+}$, is given by terms of the form $$B_n = q^{\sum_{k=1}^n{d(k)}} \, \prod_{i=1}^{n} i^{\lfloor {n}/{i} \rfloor}.$$ This is readily ascertained using the representation theorem, Theorem~\ref{th3}, and Lemma~\ref{lcp}.}
\end{remark}
}
\end{example}
\begin{example}\label{ex45} {\rm Let $q \in \mathbb{Z^+}$, $q \geq 2$ and consider $X=\{q^n\, : n\in \mathbb{N}\}$. Then the generalized factorials \cite{mb1} of this set are given simply by $n!_{_a} = \prod_{k=1}^n{(q^n-q^{k-1})}$, \cite{mb1}. The factorial set $X_q$ of this set $X$ defined by setting $B_1=q$ yields the set
\begin{eqnarray}\label{div4}
X_q = \{1, q, q^4, q^8, q^{15}, q^{21}, q^{33}, q^{41}, q^{56}, q^{69}, q^{87}, q^{99}, \ldots\},
\end{eqnarray}
whose $n$-th term is $B_n$$= q^{a(n)}$ by Lemma~\ref{lkp}, where $a(n) =\sigma (1)$$ + \ldots +$$\sigma (n)$ is (n-times) the average order of $\sigma (n)$,  (\cite{hw}, Section 18.3, p.239, p. 266). The average order of the arithmetic function $\sigma (n)$ is, in fact, the $a(n)$ defined here, its asymptotics appearing explicitly in (\cite{hw}, Theorem 324). Note that this sequence $a(n)$ appears in (\cite{njs}, id.A024916) and that $n!$ does not divide $B_n$ generally, so this set is not self-factorial. However, one may still define infinitely many other factorials on it as we have seen (see Theorem~\ref{thh9})}.
\end{example}
\begin{example} {\rm Let $q \geq 2$ be an integer and consider the integer sequence $X=\{q^{n^2}\, : n\in \mathbb{N}\}$. The factorial set $X_q$ of this set $X$ defined by setting $B_1=q$  gives the set
\begin{eqnarray}\label{div4v}
X_q = \{1, q, q^6, q^{16}, q^{37}, q^{63}, q^{113}, q^{163}, q^{248}, q^{339}, q^{469}, q^{591}, \ldots\},
\end{eqnarray}
where now the $n$-th term is $B_n = q^{a_2(n)}$ by Lemma~\ref{lkp}, where $a_2(n) =\sum_{k=1}^n \sigma_{2} (k)$ and $\sigma_{2} (k)$ represents the sum of the squares of the divisors of $k$ (\cite{hw}, p.239).}
\end{example}
The previous result generalizes nicely.
\begin{example}\label{exnk} {\rm Let $q \geq 2$, $k \geq 1$ be integers and consider the integer sequence $X=\{q^{n^k}\, : n\in \mathbb{N}\}$. In this case, the factorial set $X_q$ of this set $X$ defined as usual by setting $B_1=q$  gives the set whose  $n$-th  term is $B_n = q^{a_k(n)}$ by Lemma~\ref{lcp}, where $a_k(n) =\sum_{i=1}^n \sigma_{k} (i)$ and $\sigma_{k} (i)$ is the sum of the $k$-th powers of the divisors of $i$ (\cite{hw}, p.239).}
\end{example}

\subsection{Factorial sets of the set of primes}\label{s5}

In this section we find a factorial set for the set of primes that leads to a factorial function that is different from the one found in \cite{mb1} and describe a few of its properties.
\begin{example}\label{eX1} {\rm Let $I=\{p_i : i \in \mathbb{Z^+}\}$ be the set of primes. Setting $B_1=2$ we obtain the characterization of one of its factorial sets, i.e., $$X_I = \{2, 12, 120, 5040, 110880, 43243200, 1470268800, 1173274502400, \ldots\}$$ in the form, $X_I = \{B_n\}$ where (according to our construction),
\begin{equation}\label{pbn}
B_n = 2^n\,3^{\lfloor n/2 \rfloor}\,5^{\lfloor n/3 \rfloor}\cdots {p_i}^{\lfloor n/i \rfloor}\cdots {p_n}^{\lfloor n/n \rfloor}\\
= \prod_{i=1}^n {{p_i}^{\lfloor n/i \rfloor}}.
\end{equation}
First we note that for each $n$ the total number of prime factors of $B_n$ is equal to $d(1)+d(2)+\cdots + d(n)$. Next, this particular factorial set $X_1$ is actually contained within a class of numbers considered earlier by Ramanujan \cite{ra0}, namely the class of numbers of the form $ \prod_{i=1}^n {{p_i}^{a_i}}$ where $a_1 \geq a_2 \geq \ldots \geq a_n$, a class which includes the {\it highly composite numbers} (hcn) he had already defined in 1915. 

In addition, the superadditivity of the floor function and the representation of the ordinary factorial function as a product over primes (\cite{el}, Theorem 27) shows that for every positive integer $n$, $n! | B_n$, where $B_n$ is as in \eqref{pbn} (we omit the details). This now allows us to define an abstract factorial by writing $n!_{_{a}}  = B_n$. Since $X_I$ is a self-factorial set and there are no consecutive factorials we conclude from Lemma~\ref{th2} that
$$e_{_{a}} = 1+\sum_{n=1}^{\infty} 1/\{2^n\,3^{\lfloor n/2 \rfloor}\,5^{\lfloor n/3 \rfloor}\cdots {p_i}^{\lfloor n/i \rfloor}\cdots {p_n}^{\lfloor n/n \rfloor}\} \approx 1.5918741,$$
 is irrational. The semigroup property of abstract factorials (Proposition~\ref{semi}) implies that the sum of the reciprocals of any fixed integer  power of $B_n$ is irrational as well.}
\end{example}
The arithmetical nature of the generalized binomial coefficients (defined in Definition~\ref{def2}(2)) corresponding to the abstract factorial \eqref{pbn} inspired by the set of primes is to be noted. It follows by Proposition~\ref{probi} that

\begin{proposition}\label{sqf} The factorial function defined by $n!_a=$$B_n$ where $B_n$ is defined in \eqref{pbn} has the property that for every $n$ and for every $k,$ $0\leq k$$\leq n$, the generalized binomial coefficient $\binom{n}{k}_a$ is odd and square-free.
\end{proposition} 
\begin{remark}{\rm In 1980 Erd\"os and Graham \cite{eg} made the conjecture that the (ordinary) central binomial coefficient $\binom{2n}{n}$ is \underline {never} square-free for $n>4$. In 1985 S\'ark\"ozy \cite{as} proved this for all sufficiently large $n$, a result that was extended later by Sander \cite{jws}. Proposition~\ref{sqf} above implies the complementary result that the (generalized) central binomial coefficient $\binom{2n}{n}_a$ associated with the abstract factorial induced by the set of primes \eqref{pbn} is \underline{always} square free, for every $n$.}
\end{remark}

Now, observe that the first $6$ elements of our class $X_I$ (defined in Example~\ref{eX1})  are hcn; there is, however, little hope of finding many more due to the following result. 
\begin{proposition}\label{fff} The sequence defined by \eqref{pbn} contains only finitely many hcn.
\end{proposition}
\begin{remark} {\rm It is interesting to note that the first failure of the left side of \eqref{hcnq1} in the proof of this result is when $n=9$. Comparing all smaller hcn (i.e., those with $a_2 \leq 8$) with our sequence we see that there are no others (for a table of hcn see \cite{ra3} (pp.151-152)); thus the 6 found at the beginning of the sequence are the only ones. The sequence $B_n$ found here grows fairly rapidly: $B_n \geq 2^{n+1} p_1p_2\cdots p_n$ although this is by no means precise.}
\end{remark}
Actually more is true regarding Proposition~\ref{fff}. The next result shows that hcn are really elusive \ldots

\begin{proposition}\label{asoc} The integer sequences defined by taking {\underline {any}} of our factorial set(s) of the set of primes, even factorial sets of the factorial sets of the set of primes etc. contain only finitely many hcn.
\end{proposition}

\subsection{On factorials of the primes and abstract factorials.\label{s62}}
We consider here the question of whether the set of the factorials of the primes is an abstract factorial.

To be precise, define $f : \mathbb{N} \to \mathbb{Z^+}$ as follows:
\[f(n) = \left \{ \begin{array}{ll}
		1 ,	&	\mbox{if  \ \ $n = 0$}, \\
		1 ,	&	\mbox{if  \ \ $n =1 $}. \\
        p_{n-1}! ,	&	\mbox{if  \ \ $n \geq 2 $}. \\
		\end{array}
	\right. 
\]
The question we ask is whether $f$ is an abstract factorial? The answer seems far from obvious. A numerical search seems to indicate that the first few binomial coefficients are indeed integers (at least up to $n=50$). Indeed, use of the lower bound \cite{pd}
$$p_{n-1} > (n-1)\{\log (n-1) +\log\log (n-1)-1\}$$
for all $n \geq 7$ gives that $$ (n-1)\{\log (n-1) +\log\log (n-1)-1\}-n > 0$$ for all such $n$ (by elementary Calculus) so that $p_{n-1}> n$ for all $n \geq 7$. We conclude that $n! | f(n)$, for all $n$. 

Now consider the (abstract) binomial coefficients $$\binom{n+1}{k+1}_a = \dfrac{p_n!}{p_k!\,p_{n-k-1}!}$$ where we can assume, without loss of generality, that $n \geq 2$ and $ 1 \leq k \leq n-1$ (the remaining cases being disposed of by observation). Since $n > k$ we factor out $p_k!$ from the numerator thereby leaving a product of $p_n - p_k$ consecutive integers that are necessarily divisible by $(p_n-p_k)!$. Thus,  in order to prove that these abstract binomial coefficients are indeed integers it suffices to show that \begin{equation}\label{con}p_n \geq p_k + p_{n-k-1}, \quad\quad 1 \leq k \leq n-1,\end{equation} and all $n \geq 2$. 

On the other hand, over $50$ years ago Segal \cite{sls} proved that the Hardy-Littlewood conjecture \cite{hl} on the convexity of $\pi(x)$, i.e.,
$$\pi(x+y) \leq \pi(x)+\pi(y)$$ for all $x, y \geq 2$ is equivalent to the inequality \begin{equation}\label{hl} p_n \geq p_{n-k}+p_{k+1}-1\end{equation} for $1 \leq k \leq (n-1)/2$, $n \geq 3$, a conjecture that has not been settled yet. However, since $p_{k+1} > p_k +1$ and $p_{n-k} > p_{n-k-1}$ it follows that \eqref{hl} implies \eqref{con}. So, any counterexample to \eqref{con} also serves as a counterexample to the stated Hardy-Littlewood conjecture. Still, \eqref{con} may be true, i.e., $f(n)$ is an abstract factorial. However, settling \eqref{con} one way or another is beyond the scope of this work.

\subsection{Factorial sets of sets of highly composite numbers}\label{s6}

It turns out that there are hcn that are divisible by arbitrarily large (ordinary) factorials.
\begin{proposition}\label{ff} Let $m \in \mathbb{Z^+}$. Then there exists a highly composite number $N$ such that $m! | N$.
\end{proposition}
\begin{remark}{\rm  It is difficult to expect Proposition~\ref{ff} to be true for {\it all} hcn larger than $N$ as can be seen by considering the hcn $N=48$ where $4! | 48$ but $4!$ does not divide the next hcn, namely, $60$. However, the proof shows that Proposition~\ref{ff} is true for all those hcn larger than $N$ for which the largest prime $p$ (appearing in the prime factorization of $N$) both exceeds $e^m$ and appears in subsequents hcn's prime factorization. (This is, of course, not always the case: e.g., the largest prime in the prime decomposition of $27720$ is $11$ but the  largest such prime for the next hcn, namely $45360$, is $7$.)}
\end{remark}
 
{\bf Terminology:} We will denote by $H=\{h_n\}$ a collection of hcn with the property that $n! | h_n$ for each $n \in \mathbb{Z^+}$ (note that the existence of such a set is guaranteed by Proposition~\ref{ff}).

\begin{proposition}\label{ffffff} The factorial set $H_{h_1}$ of $H$ is self-factorial and for $k\geq 1$ the series of the reciprocals of various powers of these hcn, i.e.,
$$\sum_{n=1}^\infty 1/\{h_1^{\lfloor n/1\rfloor }h_2^{\lfloor n/2\rfloor}h_3^{\lfloor n/3\rfloor}\cdots h_n^{\lfloor n/n\rfloor}\}^k$$
is irrational.
\end{proposition}

\section{Final Remarks}\label{s8}

We add a few remarks about further irrationality results and integral representations of abstract factorials using as a basis, the Stieltjes moment problem. For background material we refer the reader to either Akhiezer \cite{nia} or Simon \cite{bs}.

We state the Stieltjes moment problem for completeness: Given a sequence $s_0, s_1, \ldots$ of real numbers to determine a measure $d\psi$ on $[0,\infty)$ such that for every $n \geq 0$,
$$s_n = \int_{0}^\infty x^n \, d\psi(x).$$
If there is one, to determine if and when it is unique and how it can be generated. One of the basic results in this area is a theorem of Carleman \cite{tc} which states that the Stieltjes moment problem has a unique solution (i.e., is {\it determinate}) provided the moments satisfy the criterion,
$$\sum_{n=1}^\infty {s_n}^{-1/2n} = + \infty.$$

In our case we consider those (necessarily positive) sequences $s_n$ generated by abstract factorials, $n!_{_a}$, for $n \geq 0$. The prototype here is the ordinary factorial, $n!$, where $d\psi(x) = e^{-x}\,dx$, which gives the classic relation between factorials and Euler's Gamma function,
$$n! = \int_{0}^\infty x^n \, e^{-x}\,dx.$$

Another less obvious example arises from a study of the abstract factorial in Example~\ref{ex44}. For any given $a \in \mathbb{Z^+}$, the abstract factorial $n!_{_a}=(an)!/a^n$ may be represented as 
$$\dfrac{(an)!}{a^n} = \int_{0}^\infty x^n \, d\psi(x),$$ where
$$d\psi(x) = \dfrac{\exp\{-(ax)^{1/a}\}}{(ax)^{1-1/a}}\, dx.$$
Here we note that each of these measures $d\psi(x)$ are also unique by the stated result of Carleman. 

The case of a general abstract factorial $n!_{_a}$ is much more difficult. Even though we know there exists a function $\phi$ of bounded variation and of finite total variation on $[0,\infty)$ such that $$n!_{_a} = \int_{0}^\infty x^n \, d\phi(x),$$
see, e.g., \cite{rpb}, the problem is whether this $\phi$ is unique let alone exhibiting such a function in this generality. 

Integral representations of abstract factorials in terms of a solution of a Stieltjes moment problem may be useful in the search for transcendence proofs for the various irrational numbers encountered here using the ideas buried in Hilbert's (1893) proof of the transcendence of $e$ using the Gamma function.

To get irrationality results of the type presented here it merely suffices to have at our disposal an abstract factorial, as then this factorial function will provide the definition of a self-factorial set. For example, the following sample theorems are an easy consequence of Theorem~\ref{thh} and the other results herein.
\begin{theorem}\label{thqn}  Let $q_n \in \mathbb{Z^+}$ be a given integer sequence satisfying $q_0=1$ and for every $n \geq 1$, $q_i\,q_j | q_n $ for all $i, j$, $1 \leq i , j \leq n$ with $i+j = n$. Then the series
$$\sum_{n=0}^\infty \frac{1}{n!q_n}$$ is irrational.
\end{theorem}

\begin{corollary}\label{thnn2} Let $f: \mathbb{Z^+}\times \mathbb{Z^+}\to \mathbb{Z^+}$ and let $f(\cdot, q)$ be concave for each $q \in \mathbb{Z^+}$.
Then, for any $q \in \mathbb{Z^+}$, $$\sum_{n=1}^\infty \frac{1}{q^{f(n,q)}n!}$$ is irrational.
\end{corollary}

In fact, binomial coefficients can, in some cases, be used to induce abstract factorials as well as one can gather from the following consequence of the previous theorem.
\begin{corollary}\label{thbcc} Let $q \in \mathbb{Z^+}$. Then
$$\sum_{n=1}^\infty {1}/{n!q^{\binom{n+q-1}{q}}}$$
is irrational.
\end{corollary}
(The alternate series counterpart of the preceding result is also irrational as usual.)

\begin{theorem}\label{thnn} Let $q \in \mathbb{Z^+}$. Then both
$$\sum_{n=1}^\infty \frac{1}{n!^{qn}} \quad \quad\quad \sum_{n=1}^\infty \frac{(-1)^n}{n!^{qn}}$$
are irrational.
\end{theorem}

\begin{example} Let $q_n = n!/2^{\lfloor n/2\rfloor}$, $n=0,1,2,\dots$. Then $q_n \in \mathbb{Z^+}$ for every $n$, $n!q_n$ is an abstract factorial and a straightforward calculation gives us that $q_i\,q_j | q_n $ for all $i, j$, $1 \leq i \leq j \leq n$ with $i+j = n$. Hence the series
\begin{eqnarray*}
\sum^\infty_{n=0} \frac{2^{\lfloor n/2\rfloor}}{n!^2} &= &\frac{1}{4}(1+\sqrt{2})\,\, I_o(\sqrt[4]{32})) + \frac{1}{4}(1-\sqrt{2})\,\, J_o(\sqrt[4]{32})\\
&\approx & 2.56279353\dots
\end{eqnarray*}
is irrational. (Here $I_o, J_o$ are Bessel functions of the first kind of order $0$.)
\end{example}

As a final result we show independently that
\begin{theorem}\label{pn} $\displaystyle \sum_{n=1}^\infty \frac{1}{p_n!}$ is irrational.
\end{theorem}
If our function $f$, defined earlier,  (basically the $n$-th prime factorial function) turns out being an abstract factorial this would not only lead to Theorem~\ref{pn} immediately but also generate other such irrationality results using products of the factorials of the $n$-th prime with other abstract factorials, as we have seen.

\section{Proofs}\label{s7}

\begin{proof}(Proposition~\ref{semi}) It suffices to prove this for any pair of abstract factorials. To this end, write $n!_{_a} =  n!_{_1}\cdot n!_{_2}$ where $n!_{_i}$, $i=1, 2$ are abstract factorials. Then $0!_{_a}=1$, and $n! | n!_{_a}$. Finally, observe that $\displaystyle \binom{n}{k}_{_a} = \prod_{i=1}^2 \binom{n}{k}_{_i}$ where, by hypothesis, each binomial coefficient on the right has integral values. 
\end{proof}

\begin{proof}(Proposition~\ref{binco}) Note that \eqref{bin1} implies that the generalized binomial coefficients of the factorial $n!_{_a}=n!h_n$ are integers. In addition, since $h_0=1$, $0!_{_a}=1$, the divisibility condition is clear. The converse is also clear and so is omitted.
\end{proof}

\begin{proof} (Lemma~\ref{10x}) Assuming the contrary we let $!_{_a}$ be such a factorial  and let $k \geq 2$ be an integer such that $r_k=r_{k+1}=1$. Since the binomial coefficient
$$\binom{k+2}{k}_{_a} =  \frac{(k+2)!_{_a}}{2!_{_a}\,k!_{_a}} = \frac{1}{2!_{_a}}  \in \mathbb{Z^+},$$
by Definition~\ref{def2}(2), this implies that $2!_{_a} | 1$ for such $k$. On the other hand, $2! | 2!_{_a}$ by Definition~\ref{def2}(3), so we get a contradiction.
\end{proof}

\begin{proof} (Lemma~\ref{11x}) Lemma~\ref{10x} guarantees that $r_{k-1}\neq 1$. Hence $r_{k-1} \geq 2$. Assume, if possible, that $r_{k-1}= 2$. Since $(k+1)!_{_a} = k!_{_a} = 2(k-1)!_{_a}$ and the generalized binomial coefficient
$$\binom{k+1}{k-1}_{_a} = \frac{(k+1)!_{_a}}{2!_{_a}\,(k-1)!_{_a}} = \frac{2}{2!_{_a}} $$ is a positive integer, $2!_{_a}$ must be equal to either $1$ or $2$. Hence, by hypothesis, it must be equal to $1$. But then by Definition~\ref{def2}(3) $2!$ must divide $2!_{_a}=1$, so we get a contradiction.
\end{proof}

\begin{proof}(Lemma~\ref{lem000})
The sequence of factorials $n!_{_a}$ is non-decreasing by Remark~\ref{nondec}, thus, in any case $\limsup_{k \to \infty}  {1}/{r_k}  \leq  1.$ Next, let  $k_n \in \mathbb{Z^+}$, be a given infinite sequence. There are then two possibilities: Either there is a subsequence, denoted again by $k_n$, such that ${k_n}!_{_a} = {(k_n + 1)}!_{_a}$ for infinitely many $n$, or every subsequence $k_n$ has the property that ${k_n}!_{_a} \neq {(k_n + 1)}!_{_a}$ except for finitely many $n$. In the first case we get \eqref{eq100}. In the second case, since ${k_n}!_{_a}$ divides ${(k_n + 1)}!_{_a}$ (by Definition~\ref{def2}) it follows that $${(k_n + 1)}!_{_a} \geq 2{k_n}!_{_a},$$ except for finitely many $n$ and this now implies \eqref{eq101}.

The final statement is supported by an example wherein $X$ is the set of all (ordinary) primes, and the factorial function is in the sense of \cite{mb1}. In this case, the explicit formula derived in (\cite{mb1}, p.793) for these factorials can be used to show sharpness when the indices in \eqref{eq101} are odd, since then $r_k=2$ for all such $k$. (See the proof of Corollary~\ref{cor311} below.)
\end{proof}

\begin{proof} (Proposition~\ref{excep}) To see that this is a factorial function we must show that the generalized binomial coefficients $\binom{n}{k}_{_{a}}$ are positive integers for $0 \leq k \leq n$ as the other two conditions in Definition~\ref{def2} are clear by construction. Putting aside the trivial cases where $k=0, k=n$ we may assume that $1 \leq k \leq n-1$.

To see that $\binom{n}{k}_{_{a}} \in  \mathbb{Z^+}$ for $k=1,2,\ldots,n-1$ we note that, by construction, the expression for $n!_{_a}$ necessarily contains two copies of each of the terms\,\, $k!_{_a}$ and $(n-k)!_{_a}$ for each such $k$ whenever $2k \neq n$. It follows that the stated binomial coefficients are integers whenever $2k \neq n$. On the other hand, if $2k=n$ the two copies of $k!_{_a}$ in the denominator are canceled by two of the respective four copies in the numerator (since now $(n-k)!_{_a}=k!_{_a}$). Observe that \eqref{eq100} holds by construction. 
\end{proof}

\begin{proof}(Theorem~\ref{thh9}) One need only apply the Definition of an abstract factorial and the construction of the $B_n$ of this section. The only part that needs a minor explanation is the integer nature of the generalized binomial coefficients. However, note that for fixed $k \in \mathbb{N}$, 
$$\binom{n}{r}_a = \binom{n}{r} \,\, \prod_{i=1}^{n+k-r} \frac{B_{r+i}}{B_i},$$
where $1 \leq r \leq n-1$, the other cases being trivial. Finally, the right hand side must be an integer since each ratio $B_{r+i}/B_i$  is also an integer, by construction.
\end{proof}

\begin{proof} (Theorem~\ref{th3}) Note that \eqref{bn} holds for the first few $n$ by inspection so we use an induction argument: Assume that $$B_i = \prod_{j=1}^i {b_j}^{\lfloor i/k \rfloor}$$ holds for all $i \leq n-1$. Since we require $B_iB_j | B_n$ for every $i,j$, $0 \leq i \leq j \leq n$ and $i+j = n$, we note that $B_iB_{n-i} | B_n$ for $i=0,1,\ldots, \lfloor n/2 \rfloor$. On the other hand if this last relation holds for all such $i$ then by the symmetry of the product involved we get $B_iB_j | B_n$ for every $i,j$, $0 \leq i \leq j \leq n$ and $i+j = n$.  Now, writing $B_n = {b_1}^{\alpha_1}\,{b_2}^{\alpha_2}\cdots {b_n}^{\alpha_n}$ where the $\alpha_i > 0$ by construction, we compare this with the expression for $B_iB_{n-i}$, that is
\begin{eqnarray*}
B_iB_{n-i} & = & \prod_{j=1}^{i} {b_j}^{\lfloor i/j \rfloor}\,  \prod_{j=1}^{n-i} {b_j}^{\lfloor (n-i)/j\rfloor},\\
& = & \prod_{j = 1}^{i} {b_j}^{\lfloor i/j \rfloor+\lfloor (n-i)/j\rfloor}\, \prod_{j = i+1}^{n-i} {b_j}^{\lfloor (n-i)/j\rfloor}.
\end{eqnarray*}
where $i \leq (n-1)/2$. Comparison of the first and last terms of this product with the expression for $B_n$ reveals that $\alpha_1=n$ and $\alpha_n=1$. For a given $j$, $1\leq j \leq n$ our construction and the induction hypothesis implies that $\alpha_i =1+ \lfloor (n-i)/i\rfloor = \lfloor n/i \rfloor $ since $B_iB_{n-i} | B_n$. This completes the induction argument.
\end{proof}

\begin{proof} (Proposition~\ref{probi})  Set aside the cases $k=0, n$ as trivial. Since  $B_n/B_k B_{n-k}=B_n/ B_{n-k}B_k$ we may assume without loss of generality that $k\geq n/2$ and that $n \geq 2$. Using the expression \eqref{bn} for $B_n$ we note that the left hand side of \eqref{bin01} may be rewritten in the form,
\begin{equation}\label{bnbk} \frac{B_n}{B_k B_{n-k}} = \prod^{n-k}_{j=1} {b_j}^{\lfloor n/j\rfloor - \lfloor k/j\rfloor - \lfloor (n-k)/j\rfloor}\cdot\prod^{k+1}_{j=n-k+1} {b_j}^{\lfloor n/j\rfloor - \lfloor k/j\rfloor}\cdot\prod^{n}_{j=k+2} {b_j}^{\lfloor n/j\rfloor}.\end{equation}
Now the first term in the first product must be 1 since $j=1$ and $n,k$ are integers. Next, since $\lfloor x+y\rfloor \leq \lfloor x\rfloor + \lfloor y\rfloor+ 1,$ for all $x, y \geq 0$, replacing $x$ by $x-y$ we get $$0\leq \lfloor x\rfloor - \lfloor y\rfloor\leq 1+ \lfloor x-y\rfloor, \quad x \geq y.$$ Hence those exponents corresponding to $j \geq 2$ in the first product are all non-negative and bounded above by 1. Furthermore, $ \lfloor (n-k)/j\rfloor = 0$ over the range $j=n-k+1,\dots,k+1$. Using this in the above display gives that the exponents in the second product are bounded above by 1. The exponents in the third product are bounded above by $\lfloor n/(k+2)\rfloor \leq1$, they are non-increasing, and bounded below by 1. Hence the exponents in the third product are all equal to 1. Thus the exponents in \eqref{bin01} are either $0$ or $1$.

The precise determination of the exponents in \eqref{bnbk} is not difficult. For a given $j$, whether $1 \leq j \leq n-k$ or $n-k+1 \leq j \leq k+1$, writing $n, k$ in base $j$ in the form $n=n_0+n_1j+n_2j^2+\dots$, etc.  we see that, 
\[\alpha_j = 
\left \{ \begin{array}{ll}
		0, &\mbox{if  $n_0-k_0 \geq 0$ }, \\ \\
		1,&\mbox{if  $n_0-k_0 < 0$}, \\
		\end{array}
	\right. 
\]
These results can be interpreted in terms of ``carries" across the radix point if required (see e.g., \cite{kw}). Finally the value $\alpha_j=1$ in the range $k+2 \leq j \leq n$.
\end{proof}

\begin{proof}(Corollary~\ref{thebi}) The assumptions imply that the generalized binomial coefficients of the factorial defined here are given by the left side of \eqref{bin01}.
\end{proof}

\begin{proof} (Lemma~\ref{th2}) The quantity $0!_{a}=1$ by definition, so we leave it out of the following discussion. Assume, on the contrary, that $e_{_a}$ is rational, that is, $E_{_a}\equiv e_{_a}-1$ is rational. Then $E_{_a} = a/b$, for some $a, b \in \mathbb{Z^+}$, $(a,b)=1$. In addition, 
$$E_{_a} - \sum_{m=1}^{k} \frac{1}{m!_{_a}} = \sum_{m=k+1}^{\infty} \frac{1}{m!_{_a}}.$$
Let $k \geq b$, $k\in \mathbb{Z^+}$ and define a number $\alpha_{_{k}}$ by setting
\begin{eqnarray}\label{eq2}
\alpha_{_{k}} &\equiv &k!_{_a} \left ( E_{_a} - \sum_{m=1}^{k} \frac{1}{m!_{_a}}\right )= k!_{_a} \left ( \frac{a}{b}- \sum_{m=1}^{k} \frac{1}{m!_{_a}}\right ) .
\end{eqnarray}
Since $k \geq b$ and $k!$ divides $k!_{_a}$ (by Definition~\ref{def2}(3)) it follows that $b$ divides $k!_{_a}$ (since $b$ divides $k!$ by our choice of $k$). Hence $k!_{_a}a/b \in \mathbb{Z^+}$. Next, for $1 \leq m \leq k$ we have that $k!_{_a}/m!_{_a} \in \mathbb{Z^+}$ (by Definition~\ref{def2}(2)). Thus, $\alpha_{_{k}} \in \mathbb{Z^+}$, for (any) $k \geq b$. 

\noindent Note that, 
\begin{eqnarray}\label{eq25}
\alpha_{_{k}} &=&k!_{_a} \sum_{m=k+1}^{\infty} \frac{1}{m!_{_a}} = k!_{_a} \left ( \frac{1}{(k+1)!_{_a}} + \frac{1}{(k+2)!_{_a}}+ \ldots \right ). 
\end{eqnarray}

First, we assume that $L < 1/2$. For $\varepsilon > 0$ so small that $L+\varepsilon < 1/2$, we choose $N$ sufficiently large so that for every $k \geq N$ we have ${k!_{_a}}/{(k+1)!_{_a}} < L+\varepsilon.$ Then it is easily verified that $$\frac{k!_{_a}}{(k+i)!_{_a}} < (L+\varepsilon)^i,$$ for every $i \geq 1$ and $k \geq N$. Since $L+\varepsilon < 1/2$ we see that 
$$\alpha_{_{k}} \leq (L+\varepsilon)\sum_{i=0}^{\infty} (L+\varepsilon)^i = \frac{L+\varepsilon}{1-(L+\varepsilon)} < 1,$$
and this leads to a contradiction.

\noindent{T}he case $L=1/2$ proceeds as above except that now we note that equality in \eqref{eq101} implies that for every $\varepsilon > 0$, there exists an $N$ such that for all $k \geq N$, $$ \frac{k!_{_a}}{(k+1)!_{_a}}  \leq 1/2 +\varepsilon.$$ Hence, for all $k \geq N$,
\begin{eqnarray}
\label{eq10a}
\alpha_{_{k}} \leq (1/2+\varepsilon)\sum_{i=0}^{\infty} (1/2+\varepsilon)^i = \frac{1/2+\varepsilon}{1-(1/2+\varepsilon)}.
\end{eqnarray}
We now fix some $\varepsilon < 1/6$ and a corresponding $N$. Then the right-side of \eqref{eq10a} is less than two. But for $k \geq N_0\equiv \max\{b,N\}$, $\alpha_{_{k}}$ is a positive integer. It follows that $\alpha_{_{k}}=1$. Using this in \eqref{eq25} we get that for every $k \geq N_0$,  
\begin{eqnarray}\label{eq11}
1 &=& k!_{_a} \left ( \frac{1}{(k+1)!_{_a}} + \frac{1}{(k+2)!_{_a}}+ \ldots \right ). 
\end{eqnarray}
Since the same argument gives that $\alpha_{_{k+1}}=1$, i.e.,
\begin{eqnarray}\label{eq12}
1 &=& (k+1)!_{_a} \left ( \frac{1}{(k+2)!_{_a}} + \frac{1}{(k+3)!_{_a}}+ \ldots \right ), 
\end{eqnarray}
comparing \eqref{eq11} and \eqref{eq12} we arrive at the relation $(k+1)!_{_a}=2k!_{_a}$, for every $k \geq N_0$. Iterating this we find that, under the assumption of equality in \eqref{eq101} we have $(k+i)!_{_a}=2^i\,k!_{_a},$ for each $i \geq 1$, and for all sufficiently large $k$. However, by Definition~\ref{def2}(3), $(k+i)!_{_a}=n\,k!_{_a}i!_{_a}$ for some $n_i\in \mathbb{Z^+}$. Hence, $n_i\,i!_{_a}=2^i$, for every $i$, for some integer $n_i$ depending on $i$. This, however, is impossible since, by Definition~\ref{def2}(4), $i!$ must divide $i!_{_a}$. Thus, $i!$ must also divide $2^i$ for every $i$ which is impossible. This contradiction proves the theorem.\\
\end{proof}

\begin{proof}(Corollary~\ref{cor311})  The prime factorization of this factorial function is given in its definition,  \eqref{mbp}. Replacing $n$, now assumed odd, by $n+1$, we see that the only contribution to $(n+1)!_{_X}$  comes from an additional factor of $2$, so that whenever $n$ is odd, we have for these factorials for the set of primes $X$ in \cite{mb1},
$$\frac{n!_X}{(n+1)!_X} = \frac{1}{2}. $$
It now follows that \eqref{eq101} is satisfied, with equality, and so by Lemma~\ref{th2}, $e_a$ is irrational.
\end{proof}

\begin{proof} (Lemma~\ref{th4}) Since $2! | 2!_{_a}$, $2!_{_a}$ must be even. There are now two cases: either $2!_{_a} \neq 2$ and this implies $2!_{_a}\geq 4$ (see Lemma~\ref{11x}), or $2!_{_a}=2$. 

{\bf Case 1:}  Let $2!_{_a} \neq 2$. We proceed as in the preceding Lemma~\ref{th2} up to \eqref{eq25}. Thus the assumption that $e_{_a}-1 $ is rational, $e_{_a} -1 =a/b$ implies that $\alpha_{_{k}} \in \mathbb{Z^+}$ \eqref{eq25} for any $k \geq b$. So, 
\begin{eqnarray}\label{eq9p}
\alpha_{_{k}} &=&k!_{_a} \sum_{n=k+1}^{\infty} \frac{1}{n!_{_a}} = k!_{_a} \left ( \frac{1}{(k+1)!_{_a}} + \frac{1}{(k+2)!_{_a}}+ \ldots \right ), \nonumber \\
& = & 1/r_{k} + 1/r_{k}r_{k+1} + \sum_{n=3}^{\infty} 1/ r_{k} r_{k+1} r_{k+2}\cdots r_{k+n-1},
\end{eqnarray}
Since the factorials have integral valued binomial coefficients we see that the product $r_{1}$$r_{2}\cdots$$r_{n-1} = n!_{_a}/1!_{_a}$ is a positive integer for every $n$. Hence, ${\binom{n+k}{k}}_{_a} \in \mathbb{Z^+}$ is equivalent to saying that $n!_{_a} | r_{k} r_{k+1} \cdots r_{k+n-1}$, for every $k \geq 0$ and $n \geq 1$. Since $n! | n!_{_a}$ for all $n$ by Definition~\ref{def2}(3), this means that 
\begin{equation}
\label{bnd}
n! | r_{k} r_{k+1} \cdots r_{k+n-1},
\end{equation}
for every integer $k \geq 0$, $n \geq 1$. 

By hypothesis there is an infinite sequence of equal consecutive factorials. Therefore, we can choose $k$ sufficiently large so that  $k \geq b$ and $r_{k+1}=1$. Then \eqref{eq9p} is satisfied for our $k$ with the $\alpha_{_{k}} $ there being a positive integer. With such a $k$ at our disposal, we now use Lemma~\ref{11x} which forces $r_{k} \geq 3$ (since $2!_{_a} \neq 2$). Using this information along with \eqref{bnd} in \eqref{eq9p} we get
\begin{eqnarray*}
\alpha_{_{k}} & \leq & 1/3 + 1/3 + \sum_{n=3}^{\infty} 1/ r_{k} r_{k+1} r_{k+2}\cdots r_{k+n-1},\\
& \leq & 2/3 + \sum_{n=3}^{\infty} 1/ n! \\
& \leq & 2/3 + e - 2 -1/2 \approx 0.8849...
\end{eqnarray*}
and this yields a contradiction.

{\bf Case 2:}  Let $2!_{_a} = 2$. We proceed as in Case 1 up to \eqref{eq9p} and then \eqref{bnd} without any changes. Once again, we choose $k \geq b$ and $r_{k+1}=1$. Since $2 = 2!_{_a} | r_{k} r_{k+1}$, we see that $r_{k}$ must be a multiple of two. If $r_{k}=2$, then \eqref{eq9p}-\eqref{bnd} together give the estimate $\alpha_{_{k}}  \leq  1/2 + 1/2 + e - 2 -1/2$ $\approx 1.218...$. However, since $\alpha_{_{k}}$ is a positive integer, we must have $\alpha_{_{k}}=1.$ Hence $r_{k}=2$ is impossible since the right side of \eqref{eq9p} must be greater than $1$. Thus, $r_{k} \geq 4$. Now using this estimate once again in \eqref{eq9p} we see that 
\begin{eqnarray}\label{eq9ppp}
1 = \alpha_{_{k}} &\leq & 1/4 + 1/4 + \sum_{n=3}^{\infty} 1/ r_{k} r_{k+1} r_{k+2}\cdots r_{k+n-1}, \\
&\leq & 1/2 + (e - 2 -1/2) \approx 0.718...
\end{eqnarray}
and there arises another final contradiction. Hence $e_{_a}$ is irrational.
\end{proof}

\begin{proof}(Theorem~\ref{thh}) This is an immediate consequence of  Lemma~\ref{th2}, Lemma~\ref{th4}, and the semigroup property.
\end{proof}

\begin{proof} (Theorem~\ref{thh2}) It suffices to prove this in the case of one factorial function with $s=1$ (by the semigroup property). This proof is simpler than the previous proof of Lemma~\ref{th2} in the unsigned case as it can be modeled on Fourier's proof of the equivalent result for the usual factorial. On the assumption that the series has a rational limit $a/b$, we let $k > b$ and then define the quantity $\alpha_k$ by
$$\alpha_k = \bigg | k!_{_a}\frac{a}{b} -  k!_{_a}\sum^k_{m=0} \frac{(-1)^{m}}{m!_{_a}} \bigg |.$$
Arguing as in Lemma~\ref{th2} we get that $\alpha_k \in \mathbb{Z^+}$ for all sufficiently large $k$.

Since the series is alternating and the sequence of factorials is non-decreasing it follows by the theory of alternating series that 
$$ 0 < \bigg | \frac{a}{b} - \sum^k_{m=0} \frac{(-1)^{m}}{m!_{_a}} \bigg | < \frac{1}{(k+1)!_{_a}}.$$
Combining the last two displays we obtain that for all sufficiently large $k$,
$$ 0 < \alpha_k < \frac{k!_{_a}}{(k+1)!_{_a}},$$
and this leads to an immediate contradiction if the factorials satisfy the alternative \eqref{eq101} in Lemma~\ref{lem000}. 

On the other hand, if the factorials satisfy the alternative \eqref{eq100} then $r_{k+1}=1$ for infinitely many $k$. We proceed as in the proof of Lemma~\ref{th4} above with minor changes. Thus, assuming the series has a rational limit $a/b$, with $a, b >0$, we can choose $k$ so large that $k > b$ so that 
\begin{eqnarray}\label{eq10p}
\beta_k \equiv k!_{_a} \sum_{m=k+1}^{\infty} \frac{(-1)^m}{m!_{_a}}  &=& k!_{_a} \left ( \frac{(-1)^{k+1}}{(k+1)!_{_a}} + \frac{(-1)^{k+2}}{(k+2)!_{_a}}+ \ldots \right ), \nonumber \\
& = & (-1)^{k+1}/r_{k} + (-1)^{k+2}/r_{k} r_{k+1} + \nonumber \\
&& + \sum_{n=3}^{\infty} (-1)^{k+n}/ r_{k} r_{k+1} r_{k+2}\cdots r_{k+n-1}.
\end{eqnarray}
But $\beta_k$ is an integer by our choice of $k$. If now $2!_{_a}  \neq 2$ (Case 1), $r_{k+1}=1$ implies that $r_k \geq 3$ and so
$$|\beta_k| \leq 1/3 + 1/3 + \sum_{n=3}^{\infty} 1/ r_{k} r_{k+1} r_{k+2}\cdots r_{k+n-1} \leq \cdots \leq 0.8849...$$
which gives a contradiction.

On the other hand, if $2!_{_a}=2$ (Case 2) then $r_{k+1}=1$ gives us that the first two terms in \eqref{eq10p} cancel out (regardless of the value of $r_k$). Hence
$$|\beta_k| \leq \sum_{n=3}^{\infty} 1/ r_{k} r_{k+1} r_{k+2}\cdots r_{k+n-1} \leq \sum_{n=3}^{\infty} 1/ n! < 1,$$
another contradiction. This completes the proof.
\end{proof}

\begin{proof}(Lemma~\ref{lns1}) According to Theorem~\ref{th3} any primitive set of the given factorial has elements $B_n$ of the form \eqref{bn} for an appropriate choice of $b_i$. Thus, if the given factorial has a primitive, then $n!_{_a}=B_n$ for all $n$. This is the case if and only if the $b_n$ are given recursively by \eqref{eqlns1}.

Conversely, if $\{b_i\}_{i=1}^\infty$ is a set of integers satisfying the divisibility condition \eqref{eqlns1}, then the set $X=\{b_1, b_2, \ldots\}$ is a primitive of this factorial.
\end{proof}

\begin{proof}(Theorem~\ref{vonm}) We define $b_1=1$, $b_i = e^{\Lambda(i)}$. The standard representation of the ordinary factorial as a product of primes (\cite{hw}, Theorem 416) gives us that 
$$\log n! =  \sum_{m \geq 1} \lfloor \frac{n}{p^m}\rfloor \log p =  \sum_{i = 1}^n \lfloor \frac{n}{i}\rfloor \log \Lambda(i) 
=\log \prod_{i=1}^n {b_i}^{\lfloor \frac{n}{i}\rfloor}
= \log B_n.$$
An application of Lemma~\ref{lns1} and \eqref{bn} now gives the conclusion. 
\end{proof}

\begin{proof}(Theorem~\ref{mbpr}) The $b_i$ being explicit, the proof is simply a matter of verification. Since \eqref{eqmbpr} is to be equal to \eqref{bn} it suffices to express the $b_i$ as products of various primes subject to their definition in the statement of the theorem. Clearly, $b_1=2$.

We observe that whenever $i \neq p^m(p-1)$ for any prime $p$ and any $m \geq 0$ the corresponding term $\lfloor n/i \rfloor$ cannot appear as a summand in \eqref{eqmbpr}. Consequently, we set $b_i=1$ in that case (as we don't want any contribution from such a term to \eqref{eqmbpr}). 

This leaves us with integers $i$, $ 1 \leq i \leq n$, that {\it can} be represented in the form $i = p^m(p-1)$ for some prime(s) and some $m \geq 0$ (the $m$'s depending on $p$ of course).  

We fix $i$. It is not hard to verify that for a given $m$ there can be at most one prime $p$ such that $i=p^m(p-1)$. For each such representation of $i$, there is a corresponding set of primes, say, $\pi_1, \pi_2, \ldots,\pi_s$, and corresponding exponents $m_1, m_2, \ldots, m_s \geq 0$ such that $i = \pi_j^{m_j}(\pi_j-1)$, $j=1,2,\ldots,s$. (For example, $4 = 2^2(2-1)=5^0(5-1)$, and there are no other such representations, so $\pi_1=2, \pi_2=5$, $m_1=2$ and $m_2=0$.)

We claim that $b_i = \pi_1\pi_2\cdots\pi_s$. (Recall that $i$ is fixed here.) Consider the term $b_i^{\lfloor \frac{n}{i}\rfloor}$  appearing in \eqref{bn}. Since
$$b_i^{\lfloor \frac{n}{i}\rfloor} = \prod_{j=1}^s \pi_j^{\lfloor \frac{n}{i}\rfloor} = \prod_{j=1}^s \pi_j^{\bigg\lfloor \dfrac{n}{\pi_j^{m_j}(\pi_j-1)}\bigg\rfloor},$$ and each multiplicand in the product must appear exactly once in the factorization \eqref{eqmbpr}, we must have all the terms in \eqref{eqmbpr} accounted for.

For if there is a prime say, $\pi_e$, from \eqref{eqmbpr} that is ``left out" of the resulting expression \eqref{bn}, there must be a corresponding denominator in the sum appearing in \eqref{eqmbpr} and so an integer $i$ having a representation in the form $\pi_e^{m_e}(\pi_e-1)$. But that prime $\pi_e$ must then appear in the resulting definition of the corresponding $b_i$. Thus all primes appearing in \eqref{eqmbpr} are accounted for in the expression \eqref{bn} and so the two quantities \eqref{eqmbpr} and \eqref{bn} must be equal.
\end{proof}

\begin{proof}(Proposition~\ref{dio}) First note that there are (infinitely many) integers $n$ that {\it cannot} be represented in the form \begin{equation}\label{dio2} n=p^m(p-1)
\end{equation}
where $p$ is a prime and $m \geq 0$. 

Thus, let $n$ satisfy \eqref{dio2} for some pair $p, m$ as required. Then, in any case, it is necessary that $p \leq n+1$. So, either $m=0$, or $m > 0$.

Let $m=0$. Then $n=p-1$, and conversely if $n$ is of the form $n=p-1$ we get one such representation (with $m=0$). Now assume that $n$ admits another representation in the form $n=q^m(q-1)$ where $q \neq p$ is another prime and $m >0$. We claim that $q$ is the largest prime factor of $n$. For otherwise, if $P>q$ is the largest such prime factor, then for some $r>0$,
there holds $n = P^r\alpha = q^m(q-1)$ where $\alpha \in \mathbb{Z^+}$ and $(P,\alpha)=1$. But since $(P, q)=1$, it is necessary that $P^r | (q-1)$. Since $P>q$ this is impossible, of course. We have thus shown that if \eqref{dio2} holds with a prime $q$ and $m >0$ then $q$ is the largest prime factor of $n$. This then gives us another possible representation of $n$ in the desired form, making this a total of at most two representations.

Let $n$ be {\it not} of the form {\it one less than a prime}, or equivalently, $m > 0$ in \eqref{dio2} for any representation of $n$ in this form. Fix such a representation \eqref{dio2}. Arguing as in the preceding case we deduce that $p$ must be the largest prime factor of $n$. Because of this, we conclude that there can be no other representation.

Thus, in conclusion, there are either no representations of a positive integer $n$ in the form \eqref{dio2} where $p$ is a prime and $m \geq 0$ (e.g., $n=7, 9,$ etc.), there is only one such representation (e.g., $n= 20, 24,$ etc.) or there are two such representations (e.g., $n=4, 18,$ etc.).
\end{proof}

\begin{proof} (Lemma~\ref{lkp}) In this generality this result is hard to find in the literature. The case $k=0$ can be found in (\cite{hw}, Theorem 320), while the case $k=1$ is referred to in (\cite{njs}, A024916). The general case can actually be found in either Basoco (\cite{mab}, eq.(16)) or Hermite (\cite{her}, p. 316).
\end{proof}

\begin{proof}  (Lemma~\ref{lcp}) Write down the list of all the divisors from $1$ to $n$ inclusively. For a given $i$, $1 \leq i \leq n$, there are $\lfloor n/i \rfloor$ multiples of the number $i$ less than or equal to $n$. Hence $i^{\lfloor n/i \rfloor}$ divides our cumulative product by definition of the latter. Taking the product over all integers $i$ shows that $\prod_{i=1}^{n} i^{\lfloor {n}/{i} \rfloor} | \alpha (n)$. But all the divisors of $\alpha(n)$ must also be in the list and so each must be a divisor of  $\prod_{i=1}^{n} i^{\lfloor {n}/{i} \rfloor}$, since there can be no omissions by the sieving method. The result follows.
\end{proof}

\begin{proof} (Lemma~\ref{asi}) For let $X_{b_1}=\{B_n\}$ be one of its factorial sets. By Theorem~\ref{th3} its terms are necessarily of the form $$B_n = {b_1}^{\lfloor n \rfloor}{b_2}^{\lfloor n/2 \rfloor}{b_3}^{\lfloor n/3 \rfloor}\cdots {b_n}^{\lfloor n/n \rfloor}.$$ Since $n! | b_n$ by hypothesis it follows that $n! | B_n$ as well, for all $n$, and so this set is self-factorial. If $b_1$ is replaced by any other element of $I$, then it is easy to see that $n! | B_n$ once again as all the exponents in the decomposition of $B_n$ are at least one.
\end{proof}

\begin{proof}(Corollary~\ref{cor46}) Since $X$ is self-factorial, $n! | B_n$ for all $n$. The elements ${B'}_n$ of $X_{B_1}$ are necessarily of the form 
$${B'}_n = {B_1}^{\lfloor n \rfloor}{B_2}^{\lfloor n/2 \rfloor}{B_3}^{\lfloor n/3 \rfloor}
\cdots {B_n}^{\lfloor n/n \rfloor}.$$ Hence, $n! | {B'}_n$ for all $n$, and this completes the proof.
\end{proof}

\begin{proof}(Proposition~\ref{sqf}) The square free part is clear on account of Proposition~\ref{probi} and the fact that the $b_i$ are primes in the representation \eqref{bn}. That the binomial coefficients must be odd is also clear since all powers of $2$ cancel out exactly by \eqref{pbn}.
\end{proof}

\begin{proof} (Proposition~\ref{fff}) This uses a deep result by Ramanujan \cite{ra0} on the structure of hcn. Once it is known that every hcn is of the form
 \begin{eqnarray}
\label{q} 
q \equiv 2^{a_2}3^{a_3}5^{a_5}\cdots p^{a_{p}}
\end{eqnarray} 
where $a_2\geq a_3 \geq a_5 \geq \cdots \geq a_{p}\geq 1$ [\cite{ra0}, III.6-8], he goes on to show that 
\begin{eqnarray}
\label{hcnq}
\lfloor{\frac{\log p}{\log \lambda}}\rfloor \leq a_\lambda \leq 2\,\lfloor{\frac{\log P}{\log \lambda}}\rfloor,
\end{eqnarray} 
for every prime index $\lambda$, (\cite{ra0}, III.6-10, eq.(54)), where $P$ is the first prime after $p$. Now set $\lambda=2$ in \eqref{hcnq} and use the fact that for the $n$-th term, $B_n$, the multiplicity of the prime $2$ is $n$, i.e., $a_2=n$. Since $p=p_n$ by the structure theorem for $B_n$, we have $P = p_{n+1}$. Since $p_n = {\rm O}(n\log n)$ for $n > 1$, (\cite{el}, Theorem 113), the right side of \eqref{hcnq} now shows that
\begin{eqnarray}\label{hcnq1} n \leq 2\,\lfloor{\frac{\log p_{n+1}}{\log 2}}\rfloor = {\rm O}(\log (n)) + {\rm O}(\log \log (n)),\end{eqnarray} which is impossible for infinitely many $n$. The result follows.
\end{proof}

\begin{proof} (Proposition~\ref{asoc}) Let $X=\{p_n\}$ be the set of primes. Recall that a factorial set is defined uniquely once we fix a value for $b_1$, some element of $X$. The choice $b_1=2$,\ldots,$b_n=p_n$ leads to the factorial set already discussed in Proposition~\ref{fff}. On the other hand, if $b_1\neq 2$ then $B_n$ can never be highly composite for any $n$ by the structure theorem for hcn. We now consider the factorial set $X_{2}$ of $X_1$ (itself the (main) factorial set of $X$ defined by setting $b_1=p_1=2$ and whose elements are given by \eqref{pbn}). The elements of $X_2$ are necessarily of the form 
\begin{eqnarray*}
B_{n,2} &=& {B_{1}}^n {B_{2}}^{\lfloor n/2\rfloor} {B_{3}}^{\lfloor n/3\rfloor} \cdots {B_{n}}^{\lfloor n/n\rfloor},\\
& = & {p_1}^n ({p_1}^2p_2)^{\lfloor n/2\rfloor} ({p_1}^3p_2p_3)^{\lfloor n/3\rfloor}\cdots ({p_1}^n {p_2}^{{\lfloor n/2\rfloor}} {p_3}^{\lfloor n/3\rfloor} \cdots p_n^{\lfloor n/n\rfloor})^{\lfloor n/n\rfloor},\\
& = & {p_1}^{\sum_{i=1}^n \lfloor i/1 \rfloor\, \lfloor n/i\rfloor} {p_2}^{\sum_{i=1}^n \lfloor i/2 \rfloor\, \lfloor n/i\rfloor}\cdots\, {p_n}^{\sum_{i=1}^n \lfloor i/n \rfloor\, \lfloor n/i\rfloor},\\
& = & {p_1}^{\sum_{i=1}^n \sigma (i)} \cdots\, p_n,
\end{eqnarray*}
where $\sigma(i)$ is the sum of the divisors of $i$ (see Lemma~\ref{lkp}). The assumption that for some $n$, $B_{n,2}$ is a hcn leads to the estimate (see \eqref{hcnq})  
\begin{eqnarray}
\label{sig}
\lfloor \log p_{n}/\log 2 \rfloor \leq \sum_{i=1}^n \sigma (i) \leq 2\lfloor \log p_{n+1}/\log 2 \rfloor.
\end{eqnarray}
However, by Theorem 324 in \cite{hw},  $\sum_{i=1}^n \sigma (i) = n^2\,\pi^2/12 + {\rm O}(n\log n).$ On the other hand, the right side of \eqref{sig} is ${\rm O}(\log n)+{\rm O}(\log \log n)$. It follows that the right hand inequality in \eqref{sig} cannot hold for infinitely many $n$, hence there can only be finitely many hcn in $X_2$. 

Observe that the more iterations we make on the factorial sets $X_1, X_2, \ldots, X_{k}$, the higher the order of the multiplicity of the prime $2$ in the decomposition of the respective terms $B_{n,k}$, and this estimate cannot be compensated by the right side of an equation of the form \eqref{sig}.
\end{proof}

\begin{proof} (Proposition~\ref{ff}) Since each prime must appear in the prime factorization of an hcn (when written as an increasing sequence) there exists a hcn of the form $$N=2^{a_2}3^{a_3}5^{a_5}\cdots p^{a_p}$$ 
with $p > e^m$ ($e = 2.718...$). Using the representation of the factorials as a product over primes we observe that $$m! | N \quad \Longleftrightarrow \quad a_{\lambda} \geq \sum_{j \geq 1} \lfloor m/\lambda^j \rfloor,$$
for every $\lambda$, where $\lambda = 2,3,5, \ldots,p$. In order to prove the latter we note that \eqref{hcnq} implies that it is sufficient to demonstrate that $$\lfloor{\frac{\log p}{\log \lambda}}\rfloor \geq   \sum_{j \geq 1} \lfloor m/\lambda^j \rfloor,$$ or since $p > e^m$ by hypothesis, that it is sufficient to show that $$\lfloor{\frac{m}{\log \lambda}}\rfloor \geq   \sum_{j \geq 1} \lfloor m/\lambda^j \rfloor,$$ for every prime $\lambda =2,3\ldots,p$. The latter, however is true on account of the estimates 
$$\lfloor{\frac{m}{\log \lambda}}\rfloor \geq \frac{m}{\lambda - 1} =  \sum_{j \geq 1} m/\lambda^j \geq \sum_{j \geq 1} \lfloor m/\lambda^j \rfloor,$$ valid for all primes $\lambda = 2,3,\dots,p$. This completes the proof.
\end{proof}

\begin{proof} (Proposition~\ref{ffffff}) Fix a factorial set $H_1 =\{h_n\}$. Then $H_1$ contains terms of the form $B_n = \prod_{j=1}^n {h_j}^{\lfloor n/j\rfloor}$ by construction where the $h_i$ are hcn in $H$. Since $n! | h_n$  Lemma~\ref{asi} implies that the factorial set $H_1$ is self-factorial. The conclusion about the irrationality now follows by Theorem~\ref{thh} since $n!_{_a} = B_n$ defines a factorial function by construction of the respective factorial sets.
\end{proof}

\begin{proof}(Theorem~\ref{thqn}) The assumptions imply that $n!_{_a}=n!q_n$ is an abstract factorial so Theorem~\ref{thh} applies and the result follows.
\end{proof}

\begin{proof}(Corollary~\ref{thnn2}) Fix $q \in \mathbb{Z^+}$, $q \geq 2$.  We define $q_0=1$ and $q_n = q^{f(n,q)}$ for $n \geq 1$. We need only verify the that the generalized binomial coefficients are integers. This, however, is a consequence of the fact that, for any $i$, $1 \leq i\leq n$, and $i+j=n$,
$$\frac{q_n}{q_i\,q_j} = q^{f(n,q) - f(i,q)-f(n-i,q)},$$ along with the concavity of $f$ in its first variable. The result is now a consequence of Theorem~\ref{thqn}.
\end{proof}

\begin{proof}(Corollary~\ref{thbcc}) Fix $q \in \mathbb{Z^+}$, and define the function $f$ by $f(n,q) = {q^{\binom{n+q-1}{q}}}$. The concavity condition is equivalent to the following inequality amongst binomial coefficients: 
\begin{equation}\label{bico}\binom{n+q-1}{q}\geq \binom{k+q-1}{q}+ \binom{n-k+q-1}{q},
\end{equation}
where $1 \leq k \leq n$. We give two proofs (one analytical, and another purely combinatorial). 

The first proof is by an induction argument on $q$. Observe that the result is true for $q=2$ as is easy to see. (Without loss of generality we assume that $2k \leq n$ throughout.) Assuming \eqref{bico} true for $q=m$, we find that 
\begin{eqnarray}
\binom{k+m}{m+1} + \binom{n-k+m}{m+1} & = & \frac{k}{m+1}\,\binom{k+m}{m}+\frac{n-k}{m+1}\,\binom{n-k+m}{m+1},\nonumber\\
& = & \frac{1}{m+1}\bigg \{(k+m)\binom{k+m-1}{m} +\nonumber\\
&& (n-k+m)\binom{n-k+m-1}{m}\bigg\}\nonumber\\
&\leq & \frac{1}{m+1}\bigg \{(k+m)\binom{k+m-1}{m} +\nonumber\\
&& (n-k+m) \bigg ( \binom{n+m-1}{m}-\binom{k+m-1}{m} \bigg )\bigg\}\nonumber\\
&\leq & \frac{2k-n}{m+1}\binom{k+m-1}{m} + \nonumber \\ 
&& \frac{m+n-k}{m+1}\binom{m+n-1}{m}. \label{bico2}
\end{eqnarray}
Using the basic identity
$$ \binom{m+n-1}{m} = \frac{m+1}{n-1} \binom{m+n-1}{m+1}$$
in \eqref{bico2} and omitting the first term therein we find, upon simplification,
\begin{eqnarray*}
\binom{k+m}{m+1} + \binom{n-k+m}{m+1} & \leq & \frac{m+n-k}{m+n} \binom{m+n}{m+1} \leq \binom{m+n}{m+1},
\end{eqnarray*}
thus completing the induction argument.

Another, much simpler, combinatorial argument due to my colleague Jason Gao follows: $\binom{n+q-1}{q}$ is the number of unordered selections (allowing repetitions)
of $m$ numbers from the set $\{1,2,\ldots,n\}$. Next, 
$$\binom{k+q-1}{q} + \binom{n-k+q-1}{q}$$
 is the number of ways of selecting $m$ numbers which are either all from
$\{1,2,\ldots,k\}$ or all from $\{k+1,k+2,\ldots,n\}$. Hence, it must be the case that \eqref{bico} holds with equality holding only in degenerate cases.
\end{proof}

\begin{proof}(Theorem~\ref{thnn}) Fix $q \in \mathbb{Z^+}$. Define the function $!_{_a}$ as follows: $0!_{_a}=1$, $n!_{_a} = n!^{qn}$. Clearly, $n! | n!_{_a}$ for all $n$, while the generalized binomial coefficients
$$\binom{n}{k}_{_a} =\left (\frac{n!}{k!}\right)^{qk} \cdot \left (\frac{n!}{(n-k)!}\right)^{q(n-k)}.$$
However, both terms on the right must be integers for $1 \leq  k \leq n$ since $r! | n!$ for all $r$ between $1$ and $n$. The result follows.
\end{proof}

\begin{proof}(Theorem~\ref{pn}) We proceed as in first part of the proof of Lemma~\ref{th2}, since the corresponding value of $L=0$ here. Assume, on the contrary, that the sum $P$ of the series is rational. Then $P = a/b$, for some $a, b \in \mathbb{Z^+}$, $(a,b)=1$. In addition, 
$$P - \sum_{m=1}^{k} \frac{1}{p_{m}!} = \sum_{m=k+1}^{\infty} \frac{1}{p_{m}!}.$$

Let $k \geq b$ and note that $p_{k}> k$ for all such $k\in \mathbb{Z^+}$. We define  $\alpha_{_{k}}$ as before by setting
\begin{eqnarray}\label{eq2}
\alpha_{_{k}} &= & p_{k}! \left ( \frac{a}{b}- \sum_{m=1}^{k} \frac{1}{p_{m}!}\right ) .
\end{eqnarray}
Since $k \geq b$ and $k!$ divides $p_{k}!$ we get that $b$ divides $p_{k}!$. Hence $p_{k}!\,a/b \in \mathbb{Z^+}$. Next, for $1 \leq m \leq k$ we have that $p_{k}!/p_{m}! \in \mathbb{Z^+}$. Thus, $\alpha_{_{k}} \in \mathbb{Z^+}$, for (any) $k \geq b$. 

\noindent As before, 
\begin{eqnarray}\label{eq9}
\alpha_{_{k}} &=&p_{k}! \sum_{m=k+1}^{\infty} \frac{1}{p_{m}!} = p_{k}! \left ( \frac{1}{p_{k+1}!} + \frac{1}{p_{k+2}!}+ \ldots \right ). 
\end{eqnarray}

Since we are in the case where $L=0$ in Lemma~\ref{th2} ({\em i.e.,} $p_{k}!/p_{k+1}! \to 0$  the rest of the proof is identical, leading to the contradiction that $\alpha_k < 1$. This contradiction completes the proof. 
\end{proof}

\section{Acknowledgments}

\address{School of Mathematics and Statistics\\
Carleton University, Ottawa, Ontario, Canada, K1S\, 5B6}\\
\email{amingare@math.carleton.ca}
\end{document}